\documentclass[11pt]{amsart}
\usepackage{t1enc}
\usepackage[latin1]{inputenc}
\usepackage[english]{babel}
\usepackage{amsmath,amsthm}
\usepackage{amsfonts}
\usepackage{latexsym}
\usepackage[dvips]{graphicx}
\usepackage{float}
\usepackage[natural]{xcolor}
\usepackage{algorithm}
\usepackage{algorithmic}
\usepackage{enumerate}
\usepackage{multirow}
\usepackage{xcolor}
\usepackage[colorlinks,linkcolor=blue]{hyperref}
\usepackage{lineno}
\usepackage{setspace}
\usepackage{amssymb}
\usepackage{fullpage}
\usepackage{tikz}
\usepackage{stfloats}
\usepackage{selinput}
\usepackage{babel}
\usepackage{nicefrac}
\usepackage{lmodern}
\usetikzlibrary{decorations.pathreplacing}
\usetikzlibrary{patterns}
\usetikzlibrary{positioning} %ÎªÁËÊµÏÖÏà¶ÔÎ»ÖÃµÄÉè¶š
\usepackage{xcolor}
\usepackage{listings}
\newtheorem{theorem}{Theorem}[section]
\newtheorem{lemma}[theorem]{Lemma}
\newtheorem{claim}{Claim}[section]

\newtheorem{conjecture}{Conjecture}[section]

\newtheorem{definition}[theorem]{Definition}
\newtheorem{proposition}[theorem]{Proposition}
\newtheorem{fact}[theorem]{Fact}

\def\int{\textrm{int}}

\begin{document}

\onehalfspacing
%\linenumbers
%\baselineskip 0.56cm
\title{Transversal Hamilton cycle in hypergraph systems}
\author{Yangyang Cheng}
\address{YC, BW and GW. School of Mathematics, Shandong University, Jinan, China \\
Email: \texttt{(YC) mathsoul@mail.sdu.edu.cn, (BW) binwang@mail.sdu.edu.cn, (GW) ghwang@sdu.edu.cn}.}
\author{Jie Han}
\address{JH. School of Mathematics and Statistics, Beijing Institute of Technology, Beijing,China, Email: \texttt{han.jie@bit.edu.cn}.}
\author{Bin Wang}
\author{Guanghui Wang}
\author{Donglei Yang}
\address{DY. Data Science Institute, Shandong University, Jinan, China \\
Email: \texttt{dlyang@sdu.edu.cn}.}

\begin{abstract}
A $k$-graph system $\emph{\textbf{H}}=\{H_i\}_{i\in[m]}$ is a family of not necessarily distinct $k$-graphs on the same $n$-vertex set $V$ and a $k$-graph $H$ on $V$ is said to be $\emph{\textbf{H}}$-transversal provided that there exists an injection $\varphi: E(H)\rightarrow [m]$ such that $e\in E(H_{\varphi(e)})$ for all $e\in E(H)$.
We show that  given $k\geq3, \gamma>0$, sufficiently large $n$ and an $n$-vertex $k$-graph system $\emph{\textbf{H}}=\{H_i\}_{i\in[n]}$, if $\delta_{k-1}(H_i)\geq(1/2+\gamma)n$ for each $i\in[n]$, then there exists an $\emph{\textbf{H}}$-transversal tight Hamilton cycle.
This extends the result of R\"{o}dl, Ruci\'{n}ski and Szemer\'{e}di [\emph{Combinatorica}, 2008] on single $k$-graphs.
%proved that for every $k\geq3, \gamma>0$, there exists $n_0$ such that for any $n>n_0$,
%each $k$-graph $H$ on $n$ vertices with $\delta_{k-1}(H)\geq(1/2+\gamma)n$ contains a Hamilton cycle, which can be seen as a natural extension of Dirac's theorem to $k$-graphs.
%We extend the above result in $k$-graphs to $k$-graph systems as follows.

%We also give some other rainbow generalizations of the Hajnal--Szemer\'{e}di theorem.

\bigskip

%\noindent {\textbf{Keywords}:}
 %Rainbow Hamilton cycle
\end{abstract}

\maketitle

\section{Introduction}

\subsection{Dirac-type problems}
The study of Hamilton cycles in graphs and hypergraphs is a central area in graph theory with a rich history.
A classical theorem of Dirac \cite{Dirac} asserts that for any $n\geq 3$, every $n$-vertex graph with minimum degree at least $n/2$ contains a Hamilton cycle.

A $k$-uniform hypergraph ($k$-graph, hereafter) $H=(V,E)$ consists of a vertex set $V$ and an edge set $E$ which is a family of $k$-element subsets of $V$, i.e. $E\subseteq\binom{V}{k}$.
%Berge \cite{MR389636} defined a Hamilton cycle in an $n$-vertex hypergraph $H$ as a cyclic ordering $v_1\cdots v_n$ such that for $i\in[n]$, there exist distinct edges $e_i$ of $E$ with $\{v_i,v_{i+1}\}\subseteq e_i$.
%The degree of a vertex $v$ in the hypergraph, is the number of edges containing $v$.
%Bermond et al. \cite{MR539937} studied the existence of Berge Hamilton cycles under the degree condition.
%In many applications, the notion of Berge Hamilton cycles appears to be not strong enough.
%Katona and Kierstead \cite{MR1671170} defined another type of cycles in hypergraphs, which has been studied extensively.
For any $S\subseteq V(H)$, the \emph{degree} of $S$ in $H$, denoted by $\deg_H(S)$, is the number of edges containing $S$.
For any integer $\ell\geq 0$, define the \emph{minimum} $\ell$-\emph{degree} $\delta_{\ell}(H)$ to be $\min\{\deg_H(S): S\in\binom{V(H)}{\ell}\}$.
%We say that an $\ell$-path $P$ connects sequences $(v_1,\ldots,v_{\ell})$ and $(v_n,\ldots,v_{n-(\ell-1)})$, which is called the \emph{ends} of $P$.
A \emph{tight Hamilton cycle} $H=(V,E)$ is a cyclic ordering of $V$ such that every interval of $k$ consecutive vertices forms an edge.
Throughout the rest of this paper, we refer to tight Hamilton cycles as Hamilton cycles.

Problems that relate the minimum degree (in general, minimum $\ell$-degree in hypergraphs) to the structure of the (hyper)graphs are often referred to as Dirac-type problems.
%We concentrate on Hamilton cycles in this paper.
%It is worth mentioning that the
The problem of determining the best possible minimum $(k-1)$-degree condition forcing Hamilton cycles in $k$-graphs, was initially researched by Katona and Kierstead \cite{MR1671170}.
They proved that every $n$-vertex $k$-graph $H$ with $\delta_{k-1}(H)>(1-\frac{1}{2k})n+4-k-\frac{5}{2k}$ admits a Hamilton cycle.
They also conjectured that the bound on the minimum $(k-1)$-degree can be reduced to roughly $n/2$, which was confirmed asymptotically by R\"{o}dl, Ruci\'{n}ski and Szemer\'{e}di in \cite{3uniform,approximate}.
The same authors gave the exact version for $k=3$ in \cite{EXACT}.

\begin{theorem}[\cite{approximate,EXACT}]
\label{main1}
Let $k\geq3, \gamma>0$ and $H$ be an $n$-vertex $k$-graph, where $n$ is sufficiently large. If $\delta_{k-1}(H)\geq(1/2+\gamma)n$, then $H$ contains a Hamilton cycle.
Furthermore, when $k=3$ it is enough to have $\delta_2(H)\ge \lfloor n/2\rfloor$.
\end{theorem}

For more problems and results on Dirac-type problems, we refer the readers to \cite{allen2021resilience,Antoniuk,J2016Loose,2017Loose,2013Minimum,2013Tight,ferber2022dirac,gould2014recent,H2010Dirac,Jie2015Minimum,2015Minimum,han2020hamiltonicity,2010Loose,D2010Hamilton,kuhn2014hamilton,D2006Loose,kuhn2012survey,37,2021Apr,mcdowell2018hamilton,2019Minimum}  and the recent surveys by R\"{o}dl and Ruci\'{n}ski \cite{30}, Simonovits and Szemer\'{e}di \cite{39} and Zhao \cite{41}.

\subsection{A transversal setting}
A $k$-\emph{graph system} $\emph{\textbf{H}}=\{H_i\}_{i\in[m]}$ is a family of not necessarily distinct $k$-graphs on the same $n$-vertex set $V$.
Note that each $H_i$ can be seen as the collection of edges with color $i$, and in this sense $\emph{\textbf{H}}$ can be regarded as an edge-colored multi-$k$-graph.
For convenience, we use $[i,j], i,j\in \mathbb{Z}$, to denote the set $\{i,i+1,\ldots,j\}$.
The set $[1,n]$ is denoted by $[n]$ in short.
\begin{definition}
Let $m$ be an integer and $\textbf{H}=\{H_i\}_{i\in[m]}$ be a $k$-graph system.
Then a $k$-graph $H$ on $V$ is said to be $\textbf{H}$-transversal provided that there exists an injection $\varphi: E(H) \rightarrow[m]$ such that $e\in E(H_{\varphi(e)})$ for each $e\in E(H)$.
\end{definition}

Since $\varphi$ is an injection, it follows that all edges of $H$ are from different members of $\emph{\textbf{H}}$.
When $m=e(H)$, $\varphi$ is a bijection and thus $H$ contains exactly one edge from each $H_i$.

Recently, the study of graph systems has received much attention.
Aharoni~\cite{MR4125343} conjectured a transversal version of the Dirac theorem, that is, for $|V|=n\geq 3$ and $\emph{\textbf{G}}=\{G_i\}_{i\in[n]}$ on $V$, if $\delta(G_i)\geq n/2$ for each $i\in[n]$, then there exists a $\emph{\textbf{G}}$-transversal Hamilton cycle.
%, that is, a cycle with edge set $\{e_1,\ldots,e_n\}$ such that $e_i\in E(G_i)$ for $i\in[n]$.
This was verified asymptotically by Cheng, Wang and Zhao~\cite{2019Rainbow}, and completely by Joos and Kim~\cite{MR4171383}.
In~\cite{from}, Bradshaw, Halasz, and Stacho strengthened the result of Joos and Kim by showing that there exist exponentially many $\emph{\textbf{G}}$-transversal Hamilton cycles.
%In~\cite{MR4125343}, it is shown that there is an \emph{n}-vertex graph system $\emph{\textbf{G}}=\{G_i\}_{i\in[3]}$, each having more than $n^2/4$ edges with no rainbow triangle and they obtain that one needs to require (roughly) at least $0.2557n^2$ edges in each $G_i$ to guarantee the existence of a rainbow triangle.
Similarly, the Hamiltonicity  of bipartite graphs due to Moon and Moser \cite{Moon} has recently been extended to the transversal setting by Bradshaw in \cite{Bradshaw}.
Recently, Ferber, Han and Mao~\cite{ferber2022dirac} gave a transversal version of the Dirac theorem in random graph systems.

%For more results on rainbow Hamilton cycles in random graphs, see
Other recent results on transversal settings include works on matchings \cite{Preprint,MR4157094,2012The,2021G,New,LU,2020A,Hongliang2018ON,2020C}, factors \cite{akbari2007rainbow,2020B,2021July} and so on.
We should mention that another related setting which studies rainbow subgraphs in edge-colored (hyper)graphs with restrictions, e.g. properly-colorings~\cite{alon2017random,chen2015long,ferber2016rainbow,frieze2014rainbow,gebauer2012rainbow,gould2022almost,gyarfas2010rainbow,gyarfas2011long,keevash2007rainbow} and general $m$-bounded colorings~\cite{albert1995multicoloured,bal2016rainbow,chen2018long,MR4055023,coulson2020rainbow,dudek2017rainbow,dudek2012rainbow,hahn1986path}.
%More Dirac-type problems in (hyper)graph systems can be found in \cite{} and we refer the readers to \cite{rainbow,Howard2017A,MR3604112,MR1652837,MR3818098} for more results on rainbow structure containments in (hyper)graph systems.

%The other is the properly colored $n$-vertex $k$-graph system.

%Given that $\textbf{\emph{H}}=\{H_1,\ldots,H_m\}$ is an $n$-vertex $k$-graph system, where $n\geq3k^2m$, if each $H_i$ is properly colored with $e(H_i)>\binom{n}{k}-\binom{n-m+1}{k}$ for $i\in[m]$, then there exists a rainbow matching.

 The main goal of this paper is to extend Theorem~\ref{main1} to the transversal setting.
 For every $k\geq3, \gamma>0$, we say that an $n$-vertex $k$-graph system $\emph{\textbf{H}}=\{H_i\}_{i\in[n]}$ is a $(k,n,\gamma)$-\emph{graph system} if $\delta_{k-1}(H_i)\geq(1/2+\gamma)n$ for each $i\in[n]$.

\begin{theorem}\label{main}
For every $k\geq3, \gamma>0$ and sufficiently large $n\in \mathbb{N}$, every $(k,n,\gamma)$-graph system $\textbf{H}=\{H_i\}_{i\in[n]}$ admits an $\textbf{H}$-transversal Hamilton cycle.
\end{theorem}

\section{Notation and Proof Strategy}
\subsection{Notation}
A \emph{tight} \emph{path} $P=v_1v_2\cdots v_t$ is a $k$-graph whose vertices can be ordered in such a way that each edge consists of $k$ consecutive vertices and two consecutive edges intersect in exactly $k-1$ vertices.
The $length$ of a tight path is the number of edges in the path.
We say that $P$ \emph{connects} $(v_1,\ldots,v_{k-1})$ and $(v_{t},\ldots,v_{t-k+2})$.
In this paper, a $(k-1)$-element sequence of distinct vertices of $V$ will be referred as a $(k-1)$-$tuple$.
The $(k-1)$-tuples $(v_1,\ldots,v_{k-1})$ and $(v_{t},\ldots,v_{t-k+2})$ are called the \emph{ends} of $P$.
Tight paths are referred as paths for convenience.

Given a $k$-graph $H$ and a $k$-graph system $\emph{\textbf{H}}=\{H_i\}_{i\in[n]}$ on the same vertex set with $H$, we define $\{i:E(H_i)\cap E(H) \neq \emptyset\}$ as the \emph{color set} of $H$, denoted by $C(H)$.
We call $P=x_1\cdots x_{2k-2}$ an $\emph{\textbf{H}}$-transversal path \emph{with color pattern} $(c_1,\ldots,c_{k-1})$ if $\{x_i,\ldots, x_{i+k-1}\}\in E(H_{c_i})$ for $i\in[k-1]$.
Let $\mathcal{P}=\{P_1,\ldots,P_m\}$ be a family of vertex-disjoint paths.
If each $P_i, i\in[m],$ is an $\emph{\textbf{H}}$-transversal path and $C(P_i)\cap C(P_j)=\emptyset$ for distinct $i,j\in[m]$, then we call this family an \emph{$\textbf{H}$-transversal family of paths}.
Denote $\bigcup_{i\in[m]}V(P_i)$ by $V(\mathcal{P})$.
The \emph{size} of $\mathcal{P}$ is the number of paths in the family.

When we write $\alpha\ll\beta$, we mean that $\alpha,\beta$ are constants in $(0, 1)$, and for every $\beta$ we have chosen, there exists $\alpha_0=\alpha_0(\beta)$ such that the subsequent arguments hold for all $\alpha\leq\alpha_0$.
While multiple constants appear in a hierarchy, they are chosen from right to left.

Besides, we require the following concentration inequalities.
\begin{proposition}[Corollary 2.3, \cite{Random}]
\label{chernoff}
Suppose that $X$ has the binomial distribution and $0<a<3/2$.
Then $\Pr(|X-\mathbb{E}X|\geq a\mathbb{E}X)\leq2e^{-a^2\mathbb{E}X/3}$.
\end{proposition}

\begin{proposition}[Corollary 2.2, \cite{2017The}]
\label{concentration}
Let $\binom{[N]}{r}$ be the set of $r$-subsets of $\{1,\ldots,N\}$ and let $h :\binom{[N]}{r}\rightarrow \mathbb{R}$ be given.
Suppose that there exists $\alpha \geq0$ such that
\[
|h(A)-h(A')|\leq \alpha
\]
for any $A, A'\in\binom{[N]}{r}$ with $|A\cap A'|=r-1$.
Let $C\subseteq [N]$ be a set of size $r$ chosen uniformly at random.
Then
\begin{equation}
\mathbb{E}e^{h(C)}=\exp(\mathbb{E}h(C)+a),
\end{equation}
where $a$ is a real constant such that $0\leq a\leq\frac{\alpha^2}{8}\min\{r, N-r\}$.
Furthermore, for any real $t>0$,
\begin{equation}
\Pr(|h(C)-\mathbb{E}h(C)|\geq t)\leq 2\exp\left(-\frac{2t^2}{\min\{r, N-r\}\alpha^2}\right).
\end{equation}
\end{proposition}

%Given a $k$-graph $H$, a \emph{fractional matching} is a function $f: E(H)\to [0,1]$, subject to the requirement that $\sum_{e:v\in e}f(e)\le 1$, for every $v\in V(H)$. Furthermore, if equality holds for every $v\in V(H)$, then we call the fractional matching \emph{perfect}.
%Let $c_{k,d}$ be the minimum $d$-degree threshold for perfect fractional matchings in $k$-graphs, namely, for every $\varepsilon>0$ and sufficiently large $n\in \mathbb N$, every $n$-vertex $k$-graph $H$ with $\delta_d(H)\ge (c_{k,d}+\eps)\binom{n-d}{k-d}$ contains a perfect fractional matching.

\subsection{Proof strategy}
In this section we give an outline of the proof of Theorem~\ref{main}.
Our proof is under the framework of the absorption method, systematised by the work of R\"{o}dl, Ruci\'{n}ski and Szemer\'{e}di~\cite{3uniform,RODL2009613}, which reduces the problem of finding a spanning subgraph to building an absorption structure and an almost spanning structure.
Tailored to our problem, the idea is to build a transversal absorbing cycle and a transversal path cover.
Moreover, the transversal absorbing cycle will be able to swallow an arbitrary leftover of vertices, a leftover of colors as well as an $\emph{\textbf{H}}$-transversal family of paths so that we obtain a transversal Hamilton cycle.
This motivates us to append the color information and the connecting technique into the transversal absorption method, which is our contribution compared with the proof of Theorem~\ref{main1} in~\cite{approximate}.
%It is stated as follows.

\begin{lemma}[\textbf{Absorbing Lemma}]\label{absorbing}
Given $k\geq3, \gamma>0$, there exists $\kappa>0$ such that the following holds for sufficiently large $n\in \mathbb{N}$.
Let $\textbf{H}=\{H_i\}_{i\in[n]}$ be a $(k,n,\gamma)$-graph system on $V$.
Then there exists an $\textbf{H}$-transversal cycle $A$ with at most $\gamma n/2$ vertices such that for any $\textbf{H}$-transversal family of paths $\mathcal{P}$ and any vertex set $U$ in $V\backslash V(A)$  with $|\mathcal{P}|,|U|\leq\kappa n$, there exists an $\textbf{H}$-transversal cycle $A'$ with vertex set $V(A)\cup U\cup V(\mathcal{P})$ and $C(A)\subseteq C(A')$.
\end{lemma}

%As for the ``transversal absorbing cycle'',
We first define two versions of absorbers as follows .

\noindent
\begin{definition}
Given a $(k,n,\gamma)$-graph system $\textbf{H}$, a vertex $x$ and a color $c$, we say that a path $P$ is a transversal absorber for $(x,c)$ in an $n$-vertex $k$-graph system if the following holds:
\begin{itemize}
  \item $x\notin V(P)$;
  \item $P=x_1\cdots x_{2k-2}$ is an $\textbf{H}$-transversal path with color pattern $(c_1,\ldots,c_{k-1})$;
  \item $x_1\cdots x_{k-1}xx_k\cdots x_{2k-2}$ is an $\textbf{H}$-transversal path with color pattern $(c,c_1,\ldots,c_{k-1})$.
\end{itemize}
\end{definition}
\vspace{12pt}
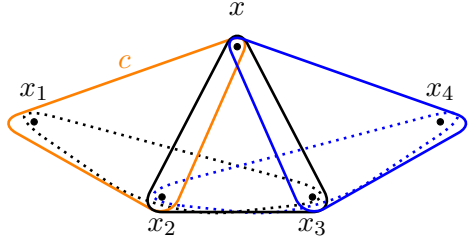
\begin{figure}[htb]
\begin{center}
\begin{tikzpicture}{center}

\filldraw [black] (-2.7,0) circle (1.2pt);
\filldraw [black](-1,-1) circle (1.2pt);
\filldraw [black] (1,-1) circle (1.2pt);
\filldraw [black] (2.7,0) circle (1.2pt);
\filldraw [black] (0,1) circle (1.2pt);
\node at (0,1.5){$x$};
\node[orange] at (-1.5,0.8){$c$};
\draw[dotted,rounded corners=0.9cm,line width =1pt] (-3.3,0.3)--(-1,-1.3)--(1.63,-1.02)--cycle;
\draw[dotted,blue,rounded corners=0.9cm,line width =1pt] (-1.63,-1.02)--(1,-1.3)--(3.4,0.3)--cycle;
\draw[orange,rounded corners=0.3cm,line width =1pt] (-3.2,0)--(-0.9,-1.3)--(0.2,1.2)--cycle;
\draw[rounded corners=0.3cm,line width =1pt] (-1.3,-1.2)--(1.3,-1.2)--(0,1.3)--cycle;
\draw[blue,rounded corners=0.3cm,line width =1pt] (0.9,-1.3)--(3.2,0)--(-0.2,1.2)--cycle;
\node at (-2.7,0.4){$x_1$};
\node at (2.7,0.4){$x_4$};
\node at (1,-1.4){$x_3$};
\node at (-1,-1.4){$x_2$};
\end{tikzpicture}
\end{center}
\caption{Absorber for $(x,c)$ when $k=3$}
\end{figure}

\vspace{12pt}

\noindent
\begin{definition}
Given a $(k,n,\gamma)$-graph system $\textbf{H}$, two disjoint $(k-1)$-tuples of vertices $\mathbf{u}=(u_1,\ldots,u_{k-1})$, $\mathbf{v}=
(v_1,\ldots,v_{k-1})$ and a $(k-1)$-tuple $(o_1,\ldots,o_{k-1})$ of colors, we say that a path $P$ is a transversal absorber for $(\mathbf{u},\mathbf{v};o_1,\ldots,o_{k-1})$ in an $n$-vertex $k$-graph system if the following holds:
\begin{itemize}
 \item $V(P)\cap \{u_1,\ldots,u_{k-1},v_1,\ldots,v_{k-1}\}=\emptyset$;
  \item $P=x_1\cdots x_{2k-2}$ is an $\textbf{H}$-transversal path with color pattern $(c_1,\ldots,c_{k-1})$;
  \item $x_1\cdots x_{k-1}u_1\cdots u_{k-1}$ and $v_1\cdots v_{k-1}x_k\cdots x_{2k-2}$ are $\textbf{H}$-transversal paths with color patterns $(c_1,\ldots,c_{k-1})$, $(o_1,\ldots,o_{k-1})$ respectively.
\end{itemize}
\end{definition}
%\begin{figure}[htb]
%  \centering
 % \includegraphics[width=6cm]{11.png}\\
 % \caption{Absorber for $(x,c)$ when $k=3$}%\label{}
%\end{figure}
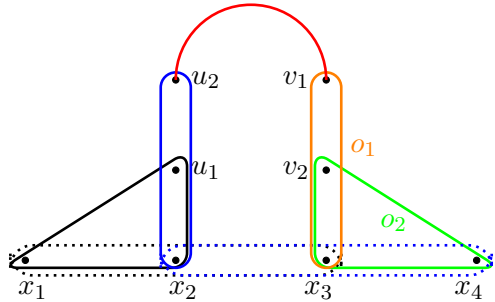
\begin{figure}[htb]
\begin{center}
\begin{tikzpicture}{center}
\filldraw [black] (-3,-1) circle (1.2pt);
\filldraw [black] (3,-1) circle (1.2pt);
\filldraw [black] (-1,0.2) circle (1.2pt);
\filldraw [black](-1,-1) circle (1.2pt);
\filldraw [black] (1,-1) circle (1.2pt);
\filldraw [black] (1,0.2) circle (1.2pt);
\filldraw [black] (-1,1.4) circle (1.2pt);
\filldraw [black] (1,1.4) circle (1.2pt);
\draw[rounded corners=0.2cm,line width =1pt] (-0.85,0.45)--(-0.85,-1.1)--(-3.3,-1.1)--cycle;
\draw[green,rounded corners=0.2cm,line width =1pt] (0.85,0.45)--(3.3,-1.1)--(0.85,-1.1)--cycle;
\draw[dotted,rounded corners=0.3cm,line width =1pt] (-3.2,-0.8) rectangle (1.2,-1.2);
\draw[dotted,blue,rounded corners=0.3cm,line width =1pt] (-1.2,-0.8) rectangle (3.2,-1.2);
\draw[blue,rounded corners=0.2cm,line width =1pt] (-1.2,1.5) rectangle (-0.8,-1.1);
\draw[orange,rounded corners=0.2cm,line width =1pt] (1.2,1.5) rectangle (0.8,-1.1);
\node at (-0.6,0.2){$u_1$};
\node at (-0.6,1.4){$u_2$};
\node at (0.6,0.2){$v_2$};
\node at (0.6,1.4){$v_1$};
\node[orange] at (1.5,0.5){$o_1$};
\node[green] at (1.9,-0.5){$o_2$};
\node at (-2.9,-1.4){$x_1$};
\node at (-0.9,-1.4){$x_2$};
\node at (0.9,-1.4){$x_3$};
\node at (2.9,-1.4){$x_4$};
\draw[red,line width =1pt](1,1.4)arc(0:180:1);
\end{tikzpicture}
\end{center}
\caption{Absorber for $((u_1,u_2),(v_1,v_2);o_1,o_2)$ when $k=3$}
\end{figure}
%\begin{figure}[htb]
%  \centering
 % % Requires \usepackage{graphicx}
  %\includegraphics[width=6cm]{22.png}\\
%  \caption{Absorber for $(\{u_1,u_2\},\{v_2,v_2\};o_1,o_2)$ when $k=3$}%\label{}
%\end{figure}

%In other words, Lemma \ref{absorbing} gives us a short rainbow cycle $A$ such that every small subset of vertices and small rainbow family of paths can be absorbed into a long rainbow cycle.
The second task is to connect the absorbers to a path.
The following lemma helps us to connect any two disjoint paths (by connecting their ends).
\begin{lemma}[\textbf{Connecting Lemma}]\label{connecting}
For every $k\geq3, \gamma>0$, there exists $c\in\mathbb{N}$ such that the following holds for sufficiently large $n\in \mathbb{N}$.
Let $\textbf{H}=\{H_i\}_{i\in [c]}$ be a $(k,n,\gamma)$-graph system and $\mathbf{u}$, $\mathbf{v}$ be two disjoint $(k-1)$-tuples of vertices.
Then, there exists an $\textbf{H}$-transversal path from $\mathbf{u}$ to $\mathbf{v}$ with at most $c+k-1$ vertices.
\end{lemma}
Given a $(k,n,\gamma)$-graph system $\emph{\textbf{H}}=\{H_i\}_{i\in [n]}$,
we need to construct an $\emph{\textbf{H}}$-transversal family of paths, covering almost all vertices of $V\backslash V(A)$ and almost all colors of $[n]\backslash C(A)$.
To achieve this, we use the regularity lemma for hypergraphs and a trick of Ferber and Kwan~\cite{ferber2022dirac}.

\begin{lemma}[\textbf{Path cover Lemma}]\label{cover}
For every $k\geq 3,$ $\gamma$, $\delta>0$, there exists $L>0$ such that the following holds for sufficiently large $n\in \mathbb{N}$. Every $(k,n,\gamma)$-graph system $\textbf{H}=\{H_i\}_{i\in [n]}$ on $V$ contains an $\textbf{H}$-transversal family of paths $\mathcal{P}$ of size at most $L$, covering at least $(1-\delta)n$ vertices of $V$.
\end{lemma}

\begin{proof}[Proof of Theorem \ref{main}]
For any $k\geq 3$ and $\gamma>0$, we choose $1/n\ll1/L,\kappa\ll\gamma,1/k$, and fix $\emph{\textbf{H}}$ to be a $(k,n,\gamma)$-graph system on $V$.
\\
\textbf{Step 1.} By Absorbing Lemma, we obtain an $\emph{\textbf{H}}$-transversal cycle $A$ with at most $\gamma n/2$ vertices such that the following property holds.
For any $\emph{\textbf{H}}$-transversal family of paths $\mathcal{P}$ and any vertex set $U$ in $V\backslash V(A)$  with $|\mathcal{P}|,|U|\leq\kappa n$,
there exists an $\emph{\textbf{H}}$-transversal cycle $A'$ with vertex set $V(A)\cup U\cup V(\mathcal{P})$ and $C(A)\subseteq C(A')$.
\\
\textbf{Step 2.} Set $\emph{\textbf{H}}'=\{H_i'\}_{i\in C'}$ where $C'=[n]\backslash C(A), H_i'=H_i[V\backslash V(A)]$  for $i\in C'$. Let $n'=n-|V(A)|$.
Note that $\emph{\textbf{H}}'$ is a $(k,n',\gamma/2)$-graph system where $n'>(1-\gamma/2)n$.
 Applying Path Cover Lemma to $\emph{\textbf{H}}'$ with $\delta=\kappa$, we obtain an $\emph{\textbf{H}}'$-transversal family of paths $\mathcal{P}=\{P_1,\ldots,P_p\}$, where $p\leq L\leq \kappa n$, covering all but at most $\kappa n'$ vertices of $V\backslash V(A)$. Denote the set of uncovered vertices by $T$. Thus, $|T|\leq\kappa n'\leq\kappa n$.
\\
\textbf{Step 3.} Using property in \textbf{Step 1}, we obtain a transversal cycle with vertex set $V(A)\cup T\cup V(\mathcal{P})$, which is actually an $\emph{\textbf{H}}$-transversal Hamilton cycle.
\end{proof}
\section{Transversal Absorption Method}
Given a vertex $x\in V$ and a color $c\in[n]$, let $\mathcal{L}(x;c)$ be the family of transversal absorbers for $(x,c)$.
Similarly, given two disjoint $(k-1)$-tuples $\textbf{u}$ and $\textbf{v}$ of $V$  and a $(k-1)$-tuple $(o_1,\ldots,o_{k-1})$ of $[n]$, let $\mathcal{L}(\textbf{u},\textbf{v};o_1,\ldots,o_{k-1})$ be the set of transversal absorbers for $(\textbf{u},\textbf{v};o_1,\ldots,o_{k-1})$. We need the following simple result.
\begin{fact}\label{neighbor}
Let $\textbf{H}=\{H_1,\ldots,H_n\}$ be a $(k,n,\gamma)$-graph system on $V$, $S$ be a $(k-1)$-subset of $V$ and $V_0\subseteq V\backslash S$. For any $i\in[n]$, we have
\[|N_{H_i}(S)\cap V_0|\geq |V_0|-\frac{1}{2}n+\gamma n+k-1.\]
In particular, for two $(k-1)$-subsets of vertices $S_1$ and $S_2$, we obtain that for any $i,j\in[n]$,
\[|N_{H_i}(S_1)\cap N_{H_j}(S_2)|\geq 2\gamma n+|S_1\cap S_2|.\]
\end{fact}
\begin{proof}
We have $|N_{H_i}(S)\cup V_0|\leq n-k+1$ and thus
\[|N_{H_i}(S)\cap V_0|\geq|V_0|+|N_{H_i}(S)|-(n+k-1)\geq|V_0|-\frac{1}{2}n+\gamma n+k-1.\]
For the second statement, we apply the first one with $S=S_1$ and $V_0=N_{H_i}(S_2)\setminus S_1$ and note that $|V_0|\geq(\frac{1}{2}+\gamma)n-(k-1-|S_1\cap S_2|)$.
\end{proof}

Next we show lower bounds on the number of absorbers in a $(k,n,\gamma)$-graph system.
\begin{proposition}\label{manyabsorber}
For any $k\geq3, \gamma>0$, there exists $\zeta>0$ such that the following holds for all sufficiently large $n\in \mathbb{N}$.
Suppose $\textbf{H}=\{H_1,\ldots,H_n\}$ is a $(k,n,\gamma)$-graph system on $V$, then $|\mathcal{L}(x;c)|\geq\zeta n^{3k-3}$ for every vertex $x\in V$ and color $c\in[n]$, $|\mathcal{L}(\mathbf{u},\mathbf{v};o_1,\ldots,o_{k-1})|\geq\zeta n^{3k-3}$
for every two disjoint $(k-1)$-tuples $\mathbf{u}$ and $\mathbf{v}$ of $V$ and a $(k-1)$-tuple $(o_1,\ldots,o_{k-1})$ of $[n]$.

\end{proposition}
\begin{proof}
Given $k,\gamma$, we choose $1/n\ll\zeta\ll\gamma/k$.
Fixing vertex $x\in V$ and color $c\in[n]$, we construct a transversal absorber $P=x_1\cdots x_{2k-2}$ with color pattern $(c_1,\ldots,c_{k-1})$ for $(x,c)$.
We choose $(c_1,\ldots,c_{k-1})$ arbitrarily, so there are $(n-1)\cdots(n-k+1)\geq2^{1-k}n^{k-1}$ choices.
Furthermore, $x_1,\ldots, x_{k-2}$ can be chosen arbitrarily in $(n-1)\cdots(n-k+2)\geq2^{2-k}n^{k-2}$ ways.
For $x_{k-1}$, there are at least $(\frac{1}{2}+\gamma)n$ choices such that $\{x_1,\ldots, x_{k-1},x\}\in E(H_c)$.
%Let $N_1:=N_{H_{c_{j-k+1}}}(\{x_{j-k+1},\ldots,x_{j-1}\})$ and $N_2:=N_{H_{c_{j-k+1}}}(\{x_{j-k+2},\ldots,x_{j-1},x\})$ for $j\in[k,2k-2]$,
By Fact \ref{neighbor},
%implies that
%\[|N_1\cap N_2|\geq 2\gamma n+k-2.
%\]
there are at least $2\gamma n+k-2$ choices for $x_j$, $j\in[k,2k-2]$, such that
$\{x_{j-k+1},\ldots,x_j\}, \{x_{j-k+2},\ldots,x_j,x\}\in E(H_{c_{j-k+1}})$.
For $j\in[k+1,2k-2]$, $x_j$ should be different from $x_1,\ldots,x_{j-k}$.
Thus,
the number of choices for each $x_j$ is at least $2\gamma n+k-2-(j-k)\geq2\gamma n$, $j\in[k,2k-2]$, yielding together at least $2^{1-k}n^{k-1}2^{2-k}n^{k-2}(\frac{1}{2}+\gamma)n(2\gamma n)^{k-1}\geq\zeta n^{3k-3}$ transversal absorbers for $(x,c)$.

Given $\textbf{u}=(u_1,\ldots,u_{k-1}),\textbf{v}=(v_1,\ldots,v_{k-1})$ and $(o_1,\ldots,o_{k-1})$, we construct a transversal absorber $P=x_1\cdots x_{2k-2}$ with color pattern $(c_1,\ldots,c_{k-1})$ for $(\textbf{u},\textbf{v};o_1,\ldots,o_{k-1})$.
There are $(n-k+1)\cdots(n-2k+2)\geq2^{1-k}n^{k-1}$ choices for $(c_1,\ldots,c_{k-1})$.
There are at least $(\frac{1}{2}+\gamma)n-(k-1)\geq\gamma n$ choices for $x_{k-1}$ such that $\{u_1,\ldots,u_{k-1},x_{k-1}\}\in E(H_{c_{k-1}})$ and $x_{k-1}$ should be different from $v_1,\ldots,v_{k-1}$.
For $x_i$, $i\in[k-2]$, there are at least $(\frac{1}{2}+\gamma)n-(2k-3)\geq\gamma n$ choices such that
 $\{u_{k-i},\ldots,u_{k-1},x_{k-1},\ldots,x_{i+1},x_i\}\in E(H_{c_i})$,
 and it should be different from $v_1,\ldots,v_{k-1},u_1,\ldots, $
 $u_{k-1-i}$.

 By Fact \ref{neighbor}, there are at least $2\gamma n$ choices for $x_k$ such that $\{x_1,\cdots,x_{k-1},x_k\}\in E(H_{c_1})$, $\{v_1,\cdots,v_{k-1},x_k\}\in E(H_{o_1})$ and it is different from $u_1,\ldots,u_{k-1}$.
 For $x_i$, $i\in[k+1,2k-2]$, the number of choices is at least $2\gamma n+k-2-(k-1+2(i-k))\geq\gamma n$, such that $\{x_{i-(k-1)},\ldots,x_i\}\in E(H_{c_{i-(k-1)}})$, $\{v_{i-(k-1)},\ldots,v_{k-1},x_k,\ldots,x_i\}\in E(H_{o_{i-(k-1)}})$ and it should be different from $u_1,\ldots,$ $u_{k-1},x_1,\ldots,x_{i-k},v_1,\ldots,v_{i-k}$.
Thus, there are at least $2^{1-k}n^{k-1}(\gamma n)^{k-1}(\gamma n)^{k-1}\geq\zeta n^{3k-3}$ transversal absorbers for $(\textbf{u},\textbf{v};o_1,\ldots,o_{k-1})$.
\end{proof}
Now we show that we can construct a family of disjoint absorbers, with all different colors.
\begin{lemma}\label{many}
For any $k\geq3$ and $\alpha,\zeta>0$, there exists $\beta>0$ such that the following holds for all sufficiently large $n\in \mathbb{N}$.
Let $\textbf{H}=\{H_1,\ldots,H_n\}$ be an $n$-vertex $k$-graph system on $V$.
If $|\mathcal{L}(x;c)|\geq\zeta n^{3k-3}$ for every vertex $x\in V$, $c\in[n]$ and $|\mathcal{L}(\mathbf{u},\mathbf{v};o_1,\ldots,o_{k-1})|\geq\zeta n^{3k-3}$ for all disjoint $(k-1)$-tuples $\mathbf{u}$ and $\mathbf{v}$ of $V$ and $(k-1)$-tuple $(o_1,\ldots,o_{k-1})$ of $[n]$, then there exists an $\textbf{H}$-transversal family $\mathcal{F}'$ of paths of length $k-1$, satisfying
\[|\mathcal{F}'|\leq\alpha n,\quad |\mathcal{F}'\cap\mathcal{L}(x;c)|\geq\beta n, \text{ and}\]
\[|\mathcal{F}'\cap\mathcal{L}(\mathbf{u},\mathbf{v};o_1,\ldots,o_{k-1})|\geq\beta n,\]
for every vertex $x\in V$, $c\in[n]$, two disjoint $(k-1)$-tuples $\mathbf{u}$ and $\mathbf{v}$ of $V$ and $(o_1,\ldots,o_{k-1})$ of $[n]$.
\end{lemma}
\begin{proof}
Let $1/n\ll\beta\ll\varepsilon\ll\alpha,\zeta$.
Each $\emph{\textbf{H}}$-transversal path $x_1x_2\cdots x_{2k-2}$ with color pattern $(c_1,\ldots,c_{k-1})$ can be viewed as a $(3k-3)$-tuple $(x_1,x_2,\ldots,x_{2k-2},c_1,\ldots,c_{k-1})$.
Choose a family $\mathcal{F}$ of $(3k-3)$-tuples from $\binom{V}{2k-2}\times\binom{[n]}{k-1}$ by including each possible $(3k-3)$-tuple independently at random with probability
\[
p=\varepsilon\frac{(n-(2k-2))!\cdot(n-(k-1))!}{(n-1)!\cdot n!}\geq\varepsilon n^{-(3k-4)}.
\]

Note that $|\mathcal{F}|$, $|\mathcal{L}(x,c)\cap\mathcal{F}|$, $|\mathcal{L}(\textbf{u},\textbf{v};o_1,\ldots,o_{k-1})\cap\mathcal{F}|$ are binomial random variables with
\[
\mathbb{E}|\mathcal{F}|=p\frac{n!\cdot n!}{(n-(2k-2))!\cdot(n-(k-1))!}=\varepsilon n,
\]
\[
\mathbb{E}|\mathcal{L}(x,c)\cap\mathcal{F}|=p|\mathcal{L}(x;c)|\geq\varepsilon\zeta n,
\]
\[
\mathbb{E}|\mathcal{L}(\textbf{u},\textbf{v};o_1,\ldots,o_{k-1})\cap\mathcal{F}|=p|\mathcal{L}(\textbf{u},\textbf{v};o_1,\ldots,o_{k-1})|\geq\varepsilon\zeta n,
\]
for every vertex $x\in V$, $c\in[n]$, two disjoint $(k-1)$-tuples $\textbf{u}$ and $\textbf{v}$ of $V$ and $(o_1,\ldots,o_{k-1})$ of $[n]$.
By Proposition \ref{chernoff}, with probability $1-o(1)$, the family $\mathcal{F}$ satisfies the following properties
\[
|\mathcal{F}|\leq2\mathbb{E}|\mathcal{F}|=2\varepsilon n\leq\alpha n,
\]
\[
|\mathcal{L}(x;c)\cap\mathcal{F}|\geq2^{-1}\mathbb{E}|\mathcal{L}(x;c)\cap\mathcal{F}|\geq2^{-1}\varepsilon\zeta n,
\]
\[
|\mathcal{L}(\textbf{u},\textbf{v};o_1,\ldots,o_{k-1})\cap\mathcal{F}|\geq2^{-1}\mathbb{E}|\mathcal{L}(\textbf{u},\textbf{v};o_1,\ldots,o_{k-1})|\geq2^{-1}\varepsilon\zeta n,
\]
for every vertex $x\in V$, $c\in[n]$, two disjoint $(k-1)$-tuples $\textbf{u}$ and $\textbf{v}$ of $V$ and $(o_1,\ldots,o_{k-1})$ of $[n]$.
We say that two $(3k-3)$-tuples $(x_1,x_2,\ldots,x_{2k-2},c_1,\ldots,c_{k-1})$ and $(y_1,y_2,\ldots,y_{2k-2},f_1,\ldots,f_{k-1})$ are \emph{intersecting} if $x_i=y_j$ for some $i,j\in[2k-2]$ or $c_m=f_{\ell}$ for some $m,\ell\in[k-1]$.
We can bound the expected number of pairs of $(3k-3)$-tuples in $\mathcal{F}$ that are intersecting from above by
\[
\frac{n!\cdot n!}{(n-(2k-2))!\cdot(n-(k-1))!}(3k-3)^2\frac{(n-1)!\cdot n!}{(n-(2k-2))!\cdot(n-(k-1))!}p^2=(3k-3)^{2}\varepsilon^2n.
\]
Thus, using Markov's inequality, we derive that with probability at least 1/2, $\mathcal{F}$ contains at most $2(3k-3)^{2}\varepsilon^2n$ intersecting pairs of $(3k-3)$-tuples.
Remove one $(3k-3)$-tuple from every intersecting pair in $\mathcal{F}$ and remove the $(3k-3)$-tuples that can not absorb any $(x,c)$ or $(\textbf{u},\textbf{v},o_1,\ldots,o_{k-1})$ where $x\in V$, $c\in[n]$, $\textbf{u}$ and $\textbf{v}$ are $(k-1)$-tuples of $V$ and $(o_1,\ldots,o_{k-1})$ is a $(k-1)$-tuple of $[n]$.
Thus the resulting subfamily, say $\mathcal{F}'$, consists of  pairwise disjoint $(3k-3)$-tuples, which satisfies
\[
|\mathcal{L}(x;c)\cap\mathcal{F}'|\geq2^{-1}\varepsilon\zeta n-2(3k-3)^{2}\varepsilon^2n\geq \beta n,
\]
for any $x\in V$, $c\in[n]$, and a similar statement holds for $|\mathcal{L}(\textbf{u},\textbf{v};o_1,\ldots,o_{k-1})\cap\mathcal{F}'|$ for any two disjoint $(k-1)$-tuples $\textbf{u}$ and $\textbf{v}$ of $V$ and a $(k-1)$-tuple $(o_1,\ldots,o_{k-1})$ of $[n]$.
Since each $(3k-3)$-tuple in $\mathcal{F}'$ induces a transversal absorber,  $\mathcal{F}'$ is an $\emph{\textbf{H}}$-transversal family of paths, where each path is of length $k-1$.
\end{proof}
Now we are ready to prove Lemma \ref{absorbing}, assuming Lemma \ref{connecting} holds.
\begin{proof}[Proof of Lemma \ref{absorbing}]
Given $1/n\ll\kappa\ll\beta\ll\alpha,\zeta\ll\gamma,1/k$, let $\emph{\textbf{H}}=\{H_i\}_{i\in[n]}$ be a $(k,n,\gamma)$-graph system on $V$.
By Proposition \ref{manyabsorber}, we obtain $|\mathcal{L}(x;c)|\geq\zeta n^{3k-3}$ for every vertex $x\in V$ and $c\in[n]$, and $|\mathcal{L}(\textbf{u},\textbf{v};o_1,\ldots,o_{k-1})|\geq\zeta n^{3k-3}$ for every two disjoint $(k-1)$-tuples $\textbf{u}$ and $\textbf{v}$ of $V$ and a $(k-1)$-tuple $(o_1,\ldots,o_{k-1})$ of $[n]$.
By Lemma \ref{many}, there is an $\emph{\textbf{H}}$-transversal family of paths $\mathcal{F}'=\{P_1,\ldots,P_{qn}\}$, where $q\leq\alpha$ and $|V(P_i)|=2k-2$ for $i\in[qn]$,
$|\mathcal{F}'\cap\mathcal{L}(x;c)|\geq\beta n$ for every vertex $x\in V$, $c\in[n]$ and
$|\mathcal{F}'\cap\mathcal{L}(\textbf{u},\textbf{v};o_1,\ldots,o_{k-1})|\geq\beta n$
 for every two disjoint $(k-1)$-tuples $\textbf{u}$ and $\textbf{v}$ of $V$ and $(o_1,\ldots,o_{k-1})$ of $[n]$.

%Denote the ends of $P_i$ by $\textbf{e}_i$ and $\textbf{f}_i$.
%Next we will connect all these paths into a cycle by iteratively connecting $\textbf{f}_i$ and $\textbf{e}_{i+1}$ by Lemma \ref{connecting} where $\textbf{e}_{q+1}=\textbf{e}_1$.

Next, we shall connect all the paths in $\mathcal{F}'$ into an $\emph{\textbf{H}}$-transversal cycle.
Suppose we have connected $P_1,\ldots,P_j$ into one path $P$, by using each time at most $\lceil8k\gamma^{-2}\rceil-(2k-2)$ vertices from outside $V(\mathcal{F}')$.
Let $\textbf{e}=(u_1,\ldots,u_{k-1})$ be an end of $P$ and $\textbf{f}=(v_1,\ldots,v_{k-1})$ be an end of $P_{j+1}$.
Let $H_i'$ be the induced subgraph of $H_i$ obtained by removing the vertices of $V(\mathcal{F}')\cup V(P)$ except $\textbf{e}$ and $\textbf{f}$.
The number of vertices removed is at most
\[
|V(\mathcal{F}')\cup V(P)|\leq (2k-2)qn+\left(\left\lceil\frac{8k}{\gamma^2}\right\rceil-(2k-2)\right)(qn-1)<\left\lceil\frac{8k}{\gamma^2}\right\rceil qn<\frac{\gamma n}{2},
\]
where the last inequality holds since $\alpha\ll\gamma$ and $k\geq 3$.

We get a $(k,n',\gamma/2)$-graph system $\emph{\textbf{H}}'=\{H_i'\}_{i\in C}$ where $C=[n]\setminus \left(C(P)\cup C(\mathcal{F}')\right)$ and $n'=|V(H_i')|$.
Taking a $(\lceil8k\gamma^{-2}\rceil-(k-1))$-subset $C'$ of $C$, we apply Lemma \ref{connecting} to $\{H_i'\}_{i\in C'}$ with $\textbf{e}'=(u_{k-1},\ldots,u_1)$ and $\textbf{f}'=(v_{k-1},\ldots,v_1)$, obtaining an $\emph{\textbf{H}}$-transversal path $P'$ connecting $\textbf{e}'$ and $\textbf{f}'$ such that $|V(P')|\leq \lceil8k\gamma^{-2}\rceil$.
Thus, $P\cup P'\cup P_{j+1}$ forms an $\emph{\textbf{H}}$-transversal path.

After connecting all paths in $\mathcal{F}'$ in a cyclic order, we obtain an $\emph{\textbf{H}}$-transversal $(k-1)$-cycle $A$ with at most
\[
(2k-2)qn+\left(\left\lceil\frac{8k}{\gamma^2}\right\rceil-(2k-2)\right)qn\leq \frac{\gamma n}{2}
\]
vertices.
Finally, fix any $\emph{\textbf{H}}$-transversal family of paths $\mathcal{P}$ and any vertex set $U$ in $V\backslash V(A)$  with $|\mathcal{P}|,|U|\leq\kappa n$.
We may assume that $U\cap V(\mathcal{P})=\emptyset$ as otherwise we just replace $U$ be $U\setminus V(\mathcal{P})$.
Since the paths in $\mathcal{P}$ are vertex disjoint, and $V(\mathcal{P})$, $U$ and $V(A)$ are pairwise disjoint, we infer that the number of colors in $\emph{\textbf{H}}$ not used in $\mathcal P$ or $A$ is at least $n-|V(A)|-\left(|V(\mathcal{P})|-(k-1)|\mathcal{P}|\right)\ge n-|V(A)|-\left(n-|V(A)|-|U|-(k-1)|\mathcal{P}|\right)= |U|+(k-1)|\mathcal P|$.
Thus, by the property of $\mathcal{F}'$ and the fact that $\kappa\ll\beta$, there is an $\emph{\textbf{H}}$-transversal cycle $A'$ with vertex set $V(A)\cup U\cup V(\mathcal{P})$ and $C(A)\subseteq C(A')$.
\end{proof}
%\section{Proof of the Connecting Lemma}
%The idea of the proof is to grow tree-like structures (called cascades) from both designated ends $e_1$ and $e_2$ until they meet, forming the desired rainbow path.
%This method can be seen in \cite{cascade,approximate}.

\section{Transversal Path cover lemma}
In this section, we prove our Path Cover Lemma, Lemma~\ref{cover}.
\subsection{Weak Regularity Lemma for Hypergraphs}
A $k$-graph $H$ is $k$-partite if there is a partition $V(H)=V_1\cup\cdots\cup V_k$ such that every edge of $H$ intersects each set $V_i$ in precisely one vertex for $i\in[k]$.
Given a $k$-partite $k$-graph $H$ on $V_1\cup\cdots\cup V_k$ and subsets $A_i\subset V_i$, $i\in[k]$, we define $e_H(A_1,\ldots,A_k)$ to be the number of edges in $H$ with one vertex in each $A_i$ and the $density$ of $H$ with respect to $(A_1,\ldots,A_k)$ as
\begin{center}
$d_H(A_1,\ldots,A_k)=\frac{e_H(A_1,\ldots,A_k)}{|A_1|\cdots|A_{k}|}$.
\end{center}

We say that a $k$-partite $k$-graph $H$ is $\varepsilon$-\emph{regular} if for all $A_i\subseteq V_i$ with $|A_i|\geq\varepsilon|V_i|$, $i\in[k]$, we have
\[
|d_H(A_1,\ldots,A_k)-d_H(V_1,\ldots,V_k)|\leq\varepsilon.
\]

We give a straightforward generalization of the graph regularity lemma.
\begin{lemma}[Weak regularity lemma for hypergraphs \cite{H2010Dirac}]
\label{regularity}
For any $k\geq2$, $\varepsilon>0$ and $t_0\in \mathbb{N}$, there exists $T_0\in \mathbb{N}$ such that the following holds.
For every $k$-graph $H$ on sufficiently large $n\in \mathbb{N}$ vertices, there is, for some $t\in \mathbb{N}$ with $t_0\leq t\leq T_0$, a partition $V(H)=V_0\cup V_1\cup\cdots\cup V_t$ such that $|V_0|\leq \varepsilon n$, $|V_1|=|V_2|=\cdots=|V_t|$
and for all but at most $\varepsilon t^k$ sets $\{i_1,\ldots,i_k\}\in\binom{[t]}{k}$, the induced $k$-partite $k$-graph $H[V_{i_1},\ldots,V_{i_k}]$ is $\varepsilon$-regular.
\end{lemma}

The partition in Lemma \ref{regularity} is called an $\varepsilon$-$regular$ $partition$ of $H$.
For an $\varepsilon$-regular partition of $H$ and $d\geq0$, we refer to the sets $V_i$, $i\in[t]$ as \emph{clusters} and define the $reduced$ $hypergraph$ $K=K(\varepsilon, d)$ with vertex set $[t]$ and $\{i_1,\ldots,i_k\}\in\binom{[t]}{k}$ being an edge if and only if $(V_{i_1},\ldots,V_{i_k})$ is $\varepsilon$-regular and $d(V_{i_1},\ldots,V_{i_k})\geq d$.
\subsection{Proof Sketch}
In this subsection we provide a proof sketch of Lemma~\ref{cover} and highlight key ideas. We need the following definition.
\begin{definition}
A hypergraph $H^*$ is a $(1,k)$-graph $((1,$k$)$-partite, in other words$)$, if there is a partition of $V(H^*)=V_1\cup V_2$ such that every edge contains exactly one vertex of $V_1$ and $k$ vertices of $V_2$.
\end{definition}

Given a partition of $V(H^*)=V_1\cup V_2$, a $(1,k-1)$-\emph{subset} $S$ of $V(H^*)$ contains one vertex in $V_1$ and $k-1$ vertices in $V_2$.
Let $\delta_{1,k-1}(H^*):=\min\{\deg_{H^*}(S): S\ {\rm is\ a\ }(1,k-1)$-${\rm subset\ of} \ V(H^*)\}$.\\
\textbf{Step 1. Construct an auxiliary $(1,k)$-graph.} Given a $(k,n,\gamma)$-graph system $\emph{\textbf{H}}=\{H_i\}_{i\in[n]}$ on $V$, we construct the auxiliary hypergraph $H^*$ with vertex set $V(H^*)=[n]\cup V$ and edge set $E(H^*)=\{\{i\}\cup e|$ $e\in E(H_i), i\in[n]\}$.
By the definition of $(k,n,\gamma)$-graph system, we have $\delta_{k-1}(H_i)\geq(1/2+\gamma)n$ for each $i\in[n]$.
Thus, $H^*$ is a $(1,k)$-graph with $\delta_{1,k-1}(H^*)\geq(1/2+\gamma)n$.\\
\textbf{Step 2. Obtain a reduced hypergraph} $K$\textbf{.}
 With an initial partition $[n]\cup V$ of $V(H^*)$, we apply the Weak Regularity Lemma (Lemma~\ref{regularity}) to $H^*$, and obtain a partition $V(H^*)=V_0^*\cup I_1\cup\cdots\cup I_{t_1}\cup W_1\cup\cdots\cup W_{t_2}$ where $|I_i|=|W_j|=m$ for every $i\in[t_1], j\in[t_2]$, $|V_0^*|\leq2\varepsilon n$.
 By moving at most $2\varepsilon n/m$ clusters to $V_0^*$ and renaming if necessary, we obtain a partition $V(H^*)=V_0\cup I_1\cup\cdots\cup I_t\cup W_1\cup\cdots\cup W_t$, where $t=\min\{t_1,t_2\}$, $|V_0|\leq 4\varepsilon n$, $I_i\subseteq [n]$ and $W_j\subseteq V$ for every $i,j\in[t]$.
Let $K$ be the reduced hypergraph for the partition with vertex set $\mathcal{I}\cup \mathcal{W}$ where $\mathcal{I}=\{I_1,\ldots, I_t\}$ and $\mathcal{W}=\{W_1,\ldots, W_t\}$.
Note that $K$ is a $(1,k)$-graph.
We will prove that $K$ almost `inherits' the $(1,k-1)$-degree condition of $H^*$ in Section \ref{tools}.\\
\textbf{Step 3. Obtain many matchings in} $K$\textbf{.} We equally split $\mathcal{I}$ into $k$ parts  $\mathcal{I}_i=\{I_{(i-1)t/k+1},\ldots,I_{it/k}\}$ for $i\in[k]$.
%Considering each $(1,k)$-graph $F_i$ of $K$ induced on $\mathcal{I}_i\cup \mathcal{W}$, we take a randomly `balanced' partition of $\mathcal{I}_i\cup \mathcal{W}$.
For each $\mathcal{I}_i\cup \mathcal{W}$, we randomly partition it to balanced smaller pieces, namely, into parts of form $\mathcal{I}'\cup\mathcal{W}'$, where $|\mathcal{I}'|=Q/k$ and $|\mathcal{W}'|=Q$.
Denote the family of vertex-disjoint $(1,k)$-subgraphs of $K$ induced on all parts from the partition of $\mathcal{I}_i\cup \mathcal{W}$ by $\mathcal{F}_i, i\in[k]$. %the collection of  vertex-disjoint $(1,k)$-subgraphs induced on all parts from the randomly balanced partition of $F_i$.
Note that the size of $\mathcal{F}_i$ is $t/Q$.
We shall see in Section \ref{random partition} that almost all members in $\mathcal{F}_i$ are `nice'  in the sense that they inherit the $(1,k-1)$-degree condition of $H^*$.
For each such member in $\mathcal{F}_i$, $i\in[k]$,  we use the following lemma, a combination of results in (\cite{2020B} Theorem 1.7) and  (\cite{MR291564} Theorem 1.2), and obtain a perfect matching.
This yields for each $i\in[k]$ a large matching (in $K$), say $M_i$, by taking the union of the resulting matchings over all members in $\mathcal{F}_i$.

\begin{lemma}[\cite{MR291564,2020B}]
\label{rainbow matching}
For every $\gamma>0,k\in \mathbb{N}$, the following holds for all sufficiently large $n\in k\mathbb{N}$.
 every $(1,k)$-graph $H$ on $[n/k]\cup V$ with $\delta_{1,k-1}(H)\geq (1/2+\gamma)n$ admits a perfect matching, where $|V|=n$.
\end{lemma}
\textbf{Step 4. Embed the paths.}
Now back to the original $(1,k)$-graph $H^*$, each matching edge in $\bigcup_{i\in[k]} M_i$ can be blown up and we obtain an $\emph{\textbf{H}}$-transversal family of paths. This is achieved in Lemma~\ref{c1}.
%Lemma~\ref{c1} enables us to `blow up' each matching $M_i,i\in[k],$ in $K$ into a partial $\emph{\textbf{H}}$-transversal family of paths.
However, note that distinct matchings $M_i, M_j$ may intersect on vertices in $\mathcal{W}$.
To overcome this, we build the
$\emph{\textbf{H}}$-transversal family of paths in $H$ piece by piece by zooming in each matching $M_i$ one by one.
 %We obtain a rainbow family of paths almost covering all vertices in Section \ref{path embedding}.

\subsection{Tools}
\label{tools}
%Recall that for a $(k,n,\gamma)$-graph system $\emph{\textbf{H}}=\{H_i\}_{i\in[n]}$ on $V$, we construct a $(1,k)$-graph $H^*$ with vertex set $[n]\cup V$ and edge set $\{i\cup e:e\in H_i, i\in[n]\}$.
%The minimum $(1,k-1)$-degree of a $(1,k)$-graph $H^*$, denoted by $\delta_{1,k-1}(H)$, is defined as the minimum of $\deg_{H^*}(S)$ over all $(1,k-1)$-subsets $S$.
%$S=\{i_0,v_1,\ldots,v_{k-1}\}$ of $V(H^*)$ where $i_0\in[n]$, $v_i\in V$ for $i\in [k-1]$.
%Note that $\delta_{1,k-1}(H^*)\geq(1/2+\gamma)n$.
 The following proposition shows that the reduced hypergraph almost inherits the minimum degree property of the original hypergraph.
\begin{proposition}\label{remaining degree}
For any $\gamma>0, k\in \mathbb{N}$, there exists $\varepsilon>0$ such that the following holds for sufficiently large $t\in \mathbb{N}$.
Given a $(1,k)$-graph $H^*$ with $\delta_{1,k-1}(H^*)\geq(1/2+\gamma)n$ and an $\varepsilon$-regular partition $V(H^*)=V_0\cup I_1\cup\cdots\cup I_t\cup W_1\cup\cdots\cup W_t$, let $K:=K(\varepsilon, \gamma/6)$ be the reduced hypergraph.
 The number of $(1,k-1)$-subsets $S$ of $V(K)$ violating $\deg_K(S)\geq(1/2+\gamma/4)t$ is at most $k\sqrt{\varepsilon} t^k$.
\end{proposition}
\begin{proof}
Let $1/t,\varepsilon\ll\gamma$.
Note that the reduced hypergraph $K(\varepsilon,\gamma/6)$ can be written as the intersection of two hypergraphs $D:=D(\gamma/6)$ and $R:=R(\varepsilon)$ both defined on the vertex set $\{I_1,\ldots, I_t, W_1,\ldots,W_t\}$ where
\begin{itemize}
  \item $D$ consists of all sets $\{I_{i_0}, W_{i_1},\ldots,W_{i_k}\}$ such that $d(I_{i_0},W_{i_1},\ldots,W_{i_k})\geq \gamma/6$,
  \item $R$ consists of all  sets $\{I_{i_0}, W_{i_1},\ldots,W_{i_k}\}$ such that $(I_{i_0},W_{i_1},\ldots,W_{i_k})$ is $\varepsilon$-regular.
\end{itemize}
For any $(1,k-1)$-set $S$, assuming $S=\{I_1, W_1,W_2,\ldots,W_{k-1}\}$, we first show that
\begin{equation}\label{13}
\deg_D(S)\geq\left(\frac{1}{2}+\frac{\gamma}{2}\right)t.
\end{equation}
Note that $n/t\geq m:=|W_i|=|I_j|$ for $i,j\in[t]$.
We now consider the number $z$ of edges in $H^*$ which intersect each of $I_{i_0}, W_{i_1},\ldots,W_{i_{k-1}}$ in exactly one vertex.
If (\ref{13}) does not hold, then from the condition on $\delta_{1,k-1}(H^*)$, we have
\begin{equation}\label{14}
tm^{k+1}\left(\frac{1}{2}+\frac{2\gamma}{3}\right)\leq m^k\left(\left(\frac{1}{2}+\gamma\right)n-(k-1)m\right)\leq z< \left(\frac{1}{2}+\frac{\gamma}{2}\right)tm^{k+1}+t\frac{\gamma}{6}m^{k+1},
\end{equation}
a contradiction.

Note that there are at most $\varepsilon t^{k+1}$ edges in $\overline{R}$ (the complement of $R$).
Let $\mathcal{S}$ be the family of all $(1,k-1)$-element subsets $S$ for which $\deg_{\overline{R}}(S)>\sqrt{\varepsilon}t$.
We have $|\mathcal{S}|\leq k\sqrt{\varepsilon}t^k$.
This, together with (\ref{13}) and $\varepsilon\ll\gamma$, implies that all but at most $k\sqrt{\varepsilon}t^k$ $(1,k-1)$-sets $S\subseteq V(K)$ satisfy $\deg_K(S)\geq\deg_D(S)-\sqrt{\varepsilon}t\geq (\frac{1}{2}+\frac{\gamma}{4})t$.
\end{proof}

%\begin{claim}\label{matching}
%If $K$ is a $(1,k)$-graph on $(1+\frac{1}{k})t$ vertices such that the number of $(1,k-1)$-element sets $S$ of vertices with $|N_K(S)|<t/k$ is only at most $\varepsilon t^k$, then there is a matching in $K$ covering at least $(1-2k\varepsilon)t$ vertices.
%\end{claim}

%\begin{lemma}\label{random}
%For any $\delta,\varepsilon>0$, given a $(1,k)$-graph $H^*$ on $[t/k]\cup V'$ with $|V'|=t$, where all but $\delta t^k/k$ of the $(1,k-1)$-tuples have degree at least $(\frac{1}{2}+\varepsilon)t$.
%Let $V_1$ be a uniformly random subset of $Q$ vertices of $V'$ and $I_1$ be a uniformly random subset of $Q/k$ vertices of $[t]$.
%Then with probability at least $1-\binom{Q}{k-1}(\delta+e^{-\Omega(\varepsilon^2Q)})$, the random induced subgraph $H^*[V_1\cup I_1]$ has minimum $(1,k-1)$-degree at least $(\frac{1}{2}+\frac{\varepsilon}{2})(Q-k+1)$.
%\end{lemma}

\subsection{Random Partition}
\label{random partition}
Ferber and Kwan \cite{ferber} showed that if we randomly partition the vertex set of a $k$-graph $H$, then the subgraph of $H$ induced on almost all parts inherits the minimum degree of $H$. Here we need such a result for our $(1,k)$-graphs, whose proof follows almost identical as that in \cite{ferber}. We include a proof for completeness.
%The result in \cite{ferber2022dirac} shows that random subgraphs of hypergraphs typically inherit minimum degree conditions.
%For our convenience, we give a stronger version.
%We states that there is a random partition of a hypergraph such that induced subgraphs on almost all parts of this partition typically inherit minimum degree conditions.
\begin{lemma}[Partition Lemma]
\label{partition}
Suppose that $k\geq3$, $\lambda,\gamma>0$, there exist $\eta>0$ and $Q\in k\mathbb{N}$ such that the following holds for $t\in Q\mathbb{N}$.
If $H^*$ is a $(1,k)$-graph on $[\frac{t}{k}]\cup V$ with $|V|=t$ where all but at most $\frac{\eta t}{k}\binom{t}{k-1}$ of the $(1,k-1)$-subsets of $V(H^*)$ have degree at least $(1/2+\gamma)(t-k+1)$, then there is a partition $V(H^*)=S_1\cup\cdots\cup S_{t/Q}$ such that all but at most $\lambda t/Q$ classes $S_i$ satisfy $\delta_{1,k-1}(H^*[S_i])\geq(1/2+\gamma/2)(Q-k+1)$ where each $S_i$ consists of a $Q/k$-subset of $[t/k]$ and a $Q$-subset of $V$.
\end{lemma}
\begin{proof}
Let $\eta\ll1/Q\ll\lambda, \gamma$.
Partition $[t/k]$ into $t/Q$ sets $I_1,\ldots,I_{t/Q}$ randomly such that $|I_i|=Q/k$ for $i\in[t/Q]$.
We randomly order $V$ as $v_1,\ldots,v_t$ and let $V_i=\{v_{(i-1)Q+1},\ldots,v_{iQ}\}$ for $i\in[t/Q]$.
Let $S_i=I_i\cup V_i$ for $i\in[t/Q]$.
Note that each $V_i$ is a random subset of $V$.
Let $M^*$ be the collection of $(1,k-1)$-subsets with degree less than $(1/2+\gamma)(t-k+1)$ in $H^*$.
Then $|M^*|\leq\frac{\eta t}{k}\binom{t}{k-1}$.
We will prove that for $i\in [t/Q]$ and every $(1,k-1)$-subset $S$ of $S_i$,
\begin{equation}
\label{equality}
\Pr\left[\deg_{H^*[S_i]}(S)<\left(\frac{1}{2}+\frac{\gamma}{2}\right)(Q-k+1)\right]\leq\eta+e^{-\Omega(\gamma^2Q)}.
\end{equation}

First note that the probability of the event $S\in M^*$, is at most $\eta$.
Now let $A_S$ be the event that $S$ is not in $M^*$.
The set $V_i\backslash S$ is equivalent to a uniformly random set of size $Q-(k-1)$ in $V\backslash S$.
Let $A$ denote the event that a vertex $v$ in $V_i\backslash S$ such that $S\cup \{v\}\in E(H^*)$.
Note that
\[\Pr[A|A_S]\geq\frac{(\frac{1}{2}+\gamma)(t-k+1)\binom{t-k}{Q-k}}{(Q-(k-1))\binom{t-(k-1)}{Q-(k-1)}}=\frac{1}{2}+\gamma,
\]
then we have
\[
\mathbb{E}\left[\deg_{H^*[S_i]}(S)|A_S\right]\geq(\frac{1}{2}+\gamma)(Q-k+1).
\]
Exchanging any element with an element outside $V_i\backslash S$ affects $\deg_{H^*[S_i]}(S)$ by at most $1$.
Fixing $i$, we apply Proposition \ref{concentration} with $S\notin M^*$, the probability that $S$ has degree less than $(1/2+\gamma/2)(Q-k+1)$ in $H^*[S_i]$ is at most
\[
2\exp\left(-\frac{2\left(\frac{\gamma}{2}(Q-k+1)\right)^2}{(Q-k+1)}\right)=e^{-\Omega(\gamma^2Q)}.
\]
Thus, (\ref{equality}) is proved.

We say that $S_i$ is $poor$ if some $(1,k-1)$-set in the induced graph $H^*[S_i]$ has  degree less than $(1/2+\gamma/2)(Q-k+1)$.
By (\ref{equality}), $\Pr[S_i \rm {\ is\ poor}]$$\leq\frac{Q}{k}\binom{Q}{k-1}(\eta+e^{-\Omega(\gamma^2Q)})$ for $i\in[t/Q]$.
Let $X$ be the number of poor classes in our partition, then $\mathbb{E}[X]\leq\frac{t}{k}\binom{Q}{k-1}(\eta+e^{-\Omega(\gamma^2Q)})$.
By Markov's inequality, we obtain
\[
\Pr\left[X\geq\lambda\frac{t}{Q}\right]\leq\frac{Q}{\lambda k}\binom{Q}{k-1}\left(\eta+e^{-\Omega(\gamma^2Q)}\right).
\]
By the choice of $\eta\ll1/Q\ll\lambda, \gamma$, it follows that
\[
\frac{Q}{k}\binom{Q}{k-1}\left(\eta+e^{-\Omega(\gamma^2Q)}\right)<\lambda,
\]
and thus $\Pr[X\geq\lambda\frac{t}{Q}]<1$.
Therefore, there is a partition $V(H^*)=S_1\cup\cdots \cup S_{t/Q}$, where $S_i=I_i\cup V_i$, such that at least $(1-\lambda)t/Q$ classes of them satisfy $\delta_{1,k-1}(H^*[S_i])\geq(1/2+\gamma/2)(Q-k+1)$.
\end{proof}
\subsection{Path Embeddings}
\label{path embedding}
Given a $(k+1)$-partite $(k+1)$-graph $H$ on $V_0\cup V_1\cup\cdots\cup V_k$, we call that a $(k-1)$-subset $S$ of $V(H)$ is \emph{legal} if $|S\cap V_i|\leq 1$ for $i\in[k]$ and $|S\cap V_0|=0$.
An \emph{expanded path} $P$ of length $t$ in $H$ is a $(k+1)$-graph with vertex set $\{c_1,\ldots,c_t\}\cup\{v_1,\ldots,v_{t+k-1}\}$ where $\{c_1,\ldots,c_t\}\subseteq V_0$, $\{v_1,\ldots,v_{t+k-1}\}\subseteq V_1\cup\cdots\cup V_k$ and edge set $\{e_1,\ldots, e_t\}$ such that $e_i= \{c_i,v_i,\ldots,v_{i+k-1}\}$.
Note that $|V(P)\cap V_j|=\lfloor\frac{t+k-1}{k}\rfloor$ or $\lceil\frac{t+k-1}{k}\rceil$ for $j\in[k]$.
%Note that a rainbow path in a $k$-graph system is a sequentially tight path in the auxiliary $(1,k)$-graph $H^*$.
\vspace{12pt}
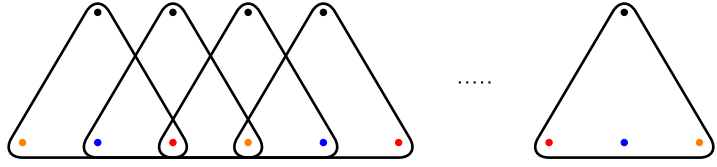
\begin{figure}[htb]
 \begin{center}
\begin{tikzpicture}
\filldraw [orange] (-4,0) circle (1.2pt);
\filldraw [blue](-3,0) circle (1.2pt);
\filldraw [red] (-2,0) circle (1.2pt);
\filldraw [orange] (-1,0) circle (1.2pt);
\filldraw [blue] (0,0) circle (1.2pt);
\filldraw [red] (1,0) circle (1.2pt);
\filldraw [orange] (5,0) circle (1.2pt);
\filldraw [blue](4,0) circle (1.2pt);
\filldraw [red] (3,0) circle (1.2pt);
\filldraw [black] (-3,1.732) circle (1.2pt);
\filldraw [black] (-2,1.732) circle (1.2pt);
\filldraw [black] (-1,1.732) circle (1.2pt);
\filldraw [black] (0,1.732) circle (1.2pt);
\filldraw [black] (4,1.732) circle (1.2pt);

\draw[dotted,line width =0.8pt] (1.8,0.8)--(2.25,0.8);

\draw[rounded corners=0.3cm,line width =1pt] (-4.3,-0.2)--(-1.7,-0.2)--(-3,2)--cycle;
\draw[rounded corners=0.3cm,line width =1pt] (-3.3,-0.2)--(-0.7,-0.2)--(-2,2)--cycle;
\draw[rounded corners=0.3cm,line width =1pt] (-2.3,-0.2)--(0.3,-0.2)--(-1,2)--cycle;
\draw[rounded corners=0.3cm,line width =1pt] (-1.3,-0.2)--(1.3,-0.2)--(0,2)--cycle;
\draw[rounded corners=0.3cm,line width =1pt] (5.3,-0.2)--(2.7,-0.2)--(4,2)--cycle;
\end{tikzpicture}
\end{center}
\caption{An expanded path in a 4-partite 4-graph (the vertices with the same color are from the same part)}
\end{figure}
\vspace{12pt}

 %If, in addition, the vertex partition satisfies

%\begin{center}
%$|V_1|\leq|V_2|\leq\cdots\leq|V_{k+1}|\leq|V_1|+1$,
%\end{center}
%then the $(k+1,k+1)$-graph will be called $equitable$.

\begin{lemma}\label{c1}
Given $c,m>0$ and $k\geq2$, every $(k+1)$-partite $(k+1)$-graph $H$ on $V_0\cup V_1\cup\cdots\cup V_k$ with at most $m$ vertices in each part and with at least $cm^{k+1}$ edges contains an expanded path of at least $cm/k$ vertices.
\end{lemma}

\begin{proof}
There are at most $k\cdot m^{k-1}$ legal $(k-1)$-subsets of $V(H)$.
We proceed the following process iteratively.
If there is a legal $(k-1)$-subset $S$, which is contained in less than $cm^2/k$ edges in the current hypergraph, then we delete all the edges containing $S$.
The process terminates at a nonempty hypergraph $H_0$ since less than $km^{k-1}(cm^2/k)=cm^{k+1}$ edges have been deleted in total.
In $H_0$, every legal $(k-1)$-subset has degree either zero or at least $cm^2/k$.

Let $P$ be a longest expanded path in $H_0$ with vertex set $\{c_1,\ldots,c_t\}\cup \{v_1,\ldots,v_{t+k-1}\}$ for some integer $t$.
We have $|V(P)\cap V_0|=t$ and $|V(P)\cap V_i|\leq t$ since each edge contains exactly one vertex of $V_i$ for each $i\in[k]$.
Consider $S_t=\{v_{t+1},\ldots,v_{t+k-1}\}$, which is a legal $(k-1)$-subset of $V(H)$.
Furthermore, $\deg_{H_0}(S_t)\geq cm^2/k$ since $S_t$ has positive degree.
All the edges containing $S_t$ must intersect $(V(P)\cap V_0)\cup(V(P)\cap V_j)$ by the maximality of $P$, where the index $j>0$ such that $S_t\cap V_j=\emptyset$.
%Thus, between any two neighbors of $e_0$ there are at least $k-1$ other vertices which also belong to $V(P)$. Hence,
Thus, we have
\begin{equation}\label{vertexnumber}
\frac{cm^2}{k}\leq |V(P)\cap V_0|\cdot|V_j|+|V(P)\cap V_j|\cdot|V_0|\leq 2tm,
\end{equation}
which implies $t\geq cm/(2k)$.
Note that $|V(P)|=t+t+k-1$ and thus $|V(P)|\geq cm/k$.
\end{proof}

The next result enables us to find a family of long vertex-disjoint expanded paths which covers almost all vertices in $V_0$ in an $\varepsilon$-regular $(k+1)$-partite $(k+1)$-graph.

\begin{lemma}\label{paths}
For any $\alpha>0,k\in \mathbb{N}$, there exists $\varepsilon>0$ such that the following holds for sufficiently large $m\in \mathbb{N}$.
Suppose $H$ is an $\varepsilon$-regular $(k+1)$-partite $(k+1)$-graph with density at least $\alpha$ and $V(H)=V_0\cup V_1\cup\cdots\cup V_k$ where $|V_0|=m$, $m/k\leq|V_i|\leq m$ for $i\in[k]$.
Then we obtain that $H$ contains a family $\mathcal{P}$ of vertex-disjoint expanded paths such that for each $P\in \mathcal{P}$, $|V(P)|\geq \varepsilon(\alpha-\varepsilon)m/k$ and $\sum_{P\in \mathcal{P}}|V(P)\cap V_0|\geq(1-2k\varepsilon)m$.
\end{lemma}
\begin{proof}
%We partition $V_i$ into $k$ classes $V_{i1},\ldots,V_{ik}$ for $i\in[k]$.
%Observe that each $(0,k-1)$-path $P$ in $H$ intersects the sets $V_1,\ldots,V_k$ almost equally, that is, for every $i,j\in[k]$, it holds that
%\[-1\leq|V(P)\cap V_i|-|V(P)\cap V_j|\leq1.\]

%Without loss of generality, we just consider the induced graph \emph{\~{H}} of $H$ on $V'=V_0\cup V_{11}\cup\cdots\cup V_{k1}$, the size of whose vertex set is $2n/(k+1)$.
Let $1/m\ll\varepsilon\ll\alpha,1/k$.
We call an expanded path $P$ \emph{good} if $|V(P)|\geq\varepsilon(\alpha-\varepsilon) m/k$.
Let $\mathcal{P}=\{P_1,\ldots,P_p\}$ be a largest family of good, vertex-disjoint expanded paths and $|V(P_i)\cap V_0|=t_i$ for $i\in[p]$. Note that $|V(P_i)\cap V_j|=\lfloor\frac{t_i+k-1}{k}\rfloor$ or $\lceil\frac{t_i+k-1}{k}\rceil$ for $i\in[p]$ and $j\in[k]$.
Suppose to the contrary that $\mathcal{P}$ covers less than $(1-2k\varepsilon) m$ vertices of $V_0$ and thus $\sum_{i\in[p]}t_i<(1-2k\varepsilon)m$.
Let $W=V(H)-\bigcup_{P\in \mathcal{P}}V(P)$ be the set of vertices uncovered by $\mathcal{P}$.
Then we have $|W\cap V_0|\geq 2k\varepsilon m$.
Hence, by the observation that $|V(P_i)\cap V_j|\leq\lceil\frac{t_i+k-1}{k}\rceil\leq\frac{t_i}{k}+2$ for each $i\in[p]$, $j\in[k]$ and the fact that $p=|\mathcal{P}|\leq (k+1)m/(\varepsilon(\alpha-\varepsilon) m/k)=k(k+1)(\varepsilon(\alpha-\varepsilon))^{-1}$, we have that
\[
|W\cap V_i|=|V_i|-|V_i\cap V(\mathcal{P})|\geq\frac{m}{k}-\sum_{i\in[p]}\left(\frac{t_i}{k}+2\right)\geq\frac{m}{k}-\frac{(1-2k\varepsilon)m}{k}-2p\geq\varepsilon m+1.
\]
Let $W_i\subseteq W\cap V_i, i\in\{0,1,\ldots,k\}$ be such that
\[
|W_0|=|W_1|=\cdots=|W_k|=\varepsilon m\geq\varepsilon|V_i|.
\]

Finally, let \emph{\^{H}} be the subhypergraph of $H$ induced on the vertex set $W_0\cup W_1\cup\cdots\cup W_k$.
%By the maximality of $\mathcal{P}$, \emph{\^{H}} is a $(k+1)$-partite $(k+1)$-graph  containing no sequentially tight path of order at least $\varepsilon(\alpha-\varepsilon)m/k$.
Since $H$ is $\varepsilon$-regular, we have
\[
d_H(W_0,W_1,\ldots,W_k)\geq d_H(V_0,V_1,\ldots,V_k)-\varepsilon\geq\alpha-\varepsilon,
\]
or equivalently,
\[
|E(\hat{H})|\geq(\alpha-\varepsilon)(\varepsilon m)^{k+1},
\]
and then Lemma \ref{c1} implies that there is an expanded path in \emph{\^{H}} on at least $\varepsilon(\alpha-\varepsilon) m/k$ vertices,  contrary to the maximality of $\mathcal{P}$.
\end{proof}
\subsection{Proof of Lemma \ref{cover}}
\begin{proof}
We choose the following parameters
\[
1/n\ll1/T_0\ll\varepsilon,1/t_0\ll1/Q\ll\lambda\ll  \delta,\gamma,1/k.
\]

Given a $(k,n,\gamma)$-graph system $\emph{\textbf{H}}=\{H_i\}_{i\in[n]}$ on $V$, we construct a $(1,k)$-graph $H^*$ with vertex set $[n]\cup V$ and edge set $\{\{i\}\cup e:e\in E(H_i), i\in[n]\}$.
With an initial partition $[n]\cup V$ of $V(H^*)$, we apply Lemma \ref{regularity} to $H^*$, and obtain a partition $V(H^*)=V_0^*\cup I_1\cup\cdots\cup I_{t_1}\cup W_1\cup\cdots\cup W_{t_2}$ where $t_0\leq t_1$, $t_2\leq T_0$, $|I_i|=|W_j|=m$ for $i\in[t_1]$ and $j\in[t_2]$, $|V_0^*|\leq2\varepsilon n$.
By moving at most $2\varepsilon n/m$ clusters to $V_0^*$ and renaming if necessary, we obtain a partition $V(H^*)=V_0\cup I_1\cup\cdots\cup I_t\cup W_1\cup\cdots\cup W_t$, where $t=\min\{t_1,t_2\}$, $|V_0|\leq 4\varepsilon n$, $I_i\subseteq [n]$ and $W_j\subseteq V$ for every $i,j\in[t]$.
Let $L=\left\lceil\frac{3kT_0}{\varepsilon(\gamma/6-\varepsilon)}\right\rceil$. %$L:=\left\lceil\frac{3kT_0}{\varepsilon(\frac{\gamma}{6}-\varepsilon)}\right\rceil$.

Let $K:=K(\varepsilon,\gamma/6)$ be the $(1,k)$-partite reduced hypergraph on $\mathcal{I}\cup \mathcal{W}$ where $\mathcal{I}=\{I_1,\ldots,I_t\}$ and $\mathcal{W}=\{W_1,\ldots,W_t\}$.
We get a family of $(1,k)$-graphs $\mathcal{F}=\{F_1,\ldots,F_{k}\}$ where $F_i=K[\{I_{(i-1)t/k+1},$ $\ldots,I_{it/k}\}\cup \mathcal{W}]$ for $i\in[k]$.

For each $i\in[k]$, applying Proposition \ref{remaining degree} and Lemma \ref{partition} to $F_i$ with $\eta=k\sqrt{\varepsilon}$, we obtain a partition $V(F_i)=S_{i,1}\cup\cdots\cup S_{i,t/Q}$ such that each $S_{i,j}$ consists of $Q/k$ vertices in $\mathcal{I}$ and $Q$ vertices in $\mathcal{W}$, and all but at most $\lambda t/Q$ classes $S_{i,j}$ satisfy $\delta_{1,k-1}(F_i[S_{i,j}])\geq(1/2+\gamma/2)(Q-k+1)$ where $j\in[t/Q]$.
We call such classes $S_{i,j}$ \emph{nice}.
Denote by $\mathcal{S}_i$ the set of indices $j$ such that $S_{i,j}$ is nice.
Then $|\mathcal{S}_i|\geq(1-\lambda)t/Q$.
%and by $\mathcal{W}_i$ the set of indices $j\in[t]$ such that $W_j$ is contained in good classes of the partition of $F_i$.
Applying Lemma \ref{rainbow matching} to each $F_i[S_{i,\ell}]$ for $i\in[k]$, $\ell\in \mathcal{S}_i$, we obtain a perfect matching $M_{i,\ell}$.
Let $M_i=\bigcup_{\ell\in\mathcal{S}_i}M_{i,\ell}$ and $M=\bigcup_{i\in[k]}M_i$.
Note that each $M_i$ is a matching in $F_i$.
For each $W_j\in \mathcal{W}$, let $p_j$ be the number of edges in $M$ that contain $W_j$, $j\in[t]$.
Then $\sum_{j\in[t]}p_j\geq k\cdot(1-\lambda) \frac{t}{Q}\cdot Q=(1-\lambda)kt$.
Next, we proceed the following process.\\
%\begin{algorithm}[H]
%\renewcommand{\algorithmicrequire}{\textbf{Input:}}
%\renewcommand\algorithmicensure {\textbf{Output:}}
%\caption{Path Embedding Process}
%\begin{algorithmic}[1]
%\REQUIRE $H^*$, $\mathcal{U}=\{I_1,\ldots,I_t\}$, $\mathcal{W}=\{W_1,\ldots,W_t\}$, $M_1,\ldots,M_k$ and $W_j^*:=W_j$ for $j\in[t]$
%\ENSURE $\mathcal{P}$
%\FOR{$i\in[k]$}
%\STATE For each $e\in M_i$, let $H_e$ be the subgraph of $H^*$ induced on the corresponding clusters constituting the edge $e$, which can be denoted by $I_e, W_{j_1(e)}^*,\ldots,W_{j_k(e)}^*$
%\STATE Applying Lemma \ref{paths} on $H_e$, we obtain a family $\mathcal{P}_e$ of vertex-disjoint $(0,k-1)$-paths that covers all but at most $2k\varepsilon m$ vertices in $I_e$
%\STATE Let $\mathcal{P}_i=\bigcup_{j\leq i}\bigcup_{e\in M_j}\mathcal{P}_e$
%\STATE Update $W_j^*$ by deleting the vertices used in $\mathcal{P}_i$ for $j\in[t]$
%\ENDFOR
%\STATE \textbf{Return} $\mathcal{P}=\mathcal{P}_k$
%\end{algorithmic}
%\end{algorithm}
\textbf{Path Embedding Process:}

Given $H^*$, $\mathcal{I}=\{I_1,\ldots,I_t\}$, $\mathcal{W}=\{W_1,\ldots,W_t\}$, $M_1,\ldots,M_k$, we initialize $W_j^*:=W_j$ for $j\in[t]$ and $i:=1$.\\
\textbf{Step 1.} For each $e\in M_i$, let $H_e$ be the subgraph of $H^*$ induced on the corresponding clusters constituting the edge $e$, we denote by $I_e, W_{j_1(e)}^*,\ldots,W_{j_k(e)}^*$ where $I_e\in\mathcal{I}$.\\
\textbf{Step 2.} Applying Lemma \ref{paths} to each $H_e$, $e\in M_i$ with $\alpha=\gamma/6$, we obtain a family $\mathcal{P}_e$ of vertex-disjoint expanded paths that covers all but at most $2k\varepsilon m$ vertices in $I_e$ and for each $P\in \mathcal{P}_e$, $|V(P)|\geq \varepsilon(\gamma/6-\varepsilon)m/k.$\\
\textbf{Step 3.} Let $\mathcal{P}_i=\bigcup_{j\leq i}\bigcup_{e\in M_j}\mathcal{P}_e$ and update $W_j^*$ by deleting the vertices used in $\mathcal{P}_i$ for $j\in[t]$. \\
\textbf{Step 4.} Update $i:=i+1$ and if $i\leq k$, go to \textbf{Step 1}; otherwise terminate the process.

After the process, we obtain $\mathcal{P}:=\mathcal{P}_k$.
It follows from the definition of $p_j$ that the size of uncovered vertices of each $W_j$ is
$$|W_j^*|=m-\sum_{W_j\in e, e\in M}|\mathcal{P}_e\cap W_j|\leq m-p_j\lfloor\frac{(1-2k\varepsilon)m+k-1}{k}\rfloor\leq m-p_j\frac{(1-2k\varepsilon)m}{k}.$$
Recall that $\sum_{j\in[t]}p_j\geq (1-\lambda)kt$. We obtain that $\mathcal{P}$ covers all but
$$|V_0|+\sum_{j\in[t]}|W_j^*|\leq4\varepsilon n+\sum_{j\in[t]}\left(m-p_j\frac{(1-2k\varepsilon)m}{k}\right)\leq\left(4(k+1)\varepsilon+\lambda\right)n\leq\delta n$$
vertices of $V$.
Moreover, since $|V(P)|\geq\frac{\varepsilon(\frac{\gamma}{6}-\varepsilon)}{k}\lfloor\frac{n}{t}\rfloor$ for each path $P\in\mathcal{P}$ and $t\leq T_0$, we have $|\mathcal{P}|<2n/(\frac{\varepsilon(\frac{\gamma}{6}-\varepsilon)}{k}\lfloor\frac{n}{t}\rfloor)<L$.
Finally, observe that $\mathcal{P}$ gives rise to an $\emph{\textbf{H}}$-transversal family of paths which completes the proof.
\end{proof}
\begin{figure}[htb]
\begin{center}
\begin{tikzpicture}[decoration=brace]
\draw[line width =0.8pt, rounded corners](-4.7,-4)rectangle(-3.3, -8);
\draw[line width =0.8pt, rounded corners](-2,-4)rectangle(0, -8);
\node at (-2.5,-8.8){$F_1$(the black dotted lines represent nice classes)};
\node at (5.6,-8.8){$F_1[S_{1,1}]$(with a perfect matching)};
\draw[line width =0.8pt,rounded corners,dotted](-4.6,-4.1)rectangle(-0.1, -4.6);
\draw[line width =0.8pt,rounded corners,dotted](-4.6,-4.7)rectangle(-0.1, -5.2);
\draw[line width =0.8pt,rounded corners,dotted](-4.6,-6.2)rectangle(-0.1, -6.7);
\draw[line width =0.8pt,rounded corners,dotted,red](-4.6,-6.8)rectangle(-0.1, -7.3);
\draw[line width =0.8pt,rounded corners,dotted,red](-4.6,-7.4)rectangle(-0.1, -7.9);
\draw[dotted,line width =0.8pt] (-2.6,-5.4)--(-2.6,-5.85);
\node at (-2.5,-4.4){$S_{1,1}$};
\node at (-2.5,-5){$S_{1,2}$};
\draw[line width =0.8pt, rounded corners](3.4,-4)rectangle(4.7, -8);
\draw[line width =0.8pt, rounded corners](6,-4)rectangle(7.9, -8);
\draw[line width =0.8pt,rounded corners,dotted](3.5,-4.1)rectangle(7.8, -4.6);
\draw[line width =0.8pt,rounded corners,dotted](3.5,-4.7)rectangle(7.8, -5.2);
\draw[line width =0.8pt,rounded corners,dotted](3.5,-6.8)rectangle(7.8, -7.3);
\draw[line width =0.8pt,rounded corners,dotted](3.5,-7.4)rectangle(7.8, -7.9);

\draw[dotted,line width =0.8pt] (5.4,-5.6)--(5.4,-6.05);
\filldraw [black] (4,-4.3) circle (1.2pt);
\filldraw [black] (6.2,-4.3) circle (1.2pt);
\filldraw [black] (7.6,-4.3) circle (1.2pt);
\draw[dotted,line width =0.8pt] (6.7,-4.3)--(7.15,-4.3);
\filldraw [black] (4,-4.9) circle (1.2pt);
\filldraw [black] (6.2,-4.9) circle (1.2pt);
\filldraw [black] (7.6,-4.9) circle (1.2pt);
\draw[dotted,line width =0.8pt] (6.7,-4.9)--(7.15,-4.9);
\filldraw [black] (4,-7) circle (1.2pt);
\filldraw [black] (6.2,-7) circle (1.2pt);
\filldraw [black] (7.6,-7) circle (1.2pt);
\draw[dotted,line width =0.8pt] (6.7,-7)--(7.15,-7);
\filldraw [black] (4,-7.65) circle (1.2pt);
\filldraw [black] (6.2,-7.65) circle (1.2pt);
\filldraw [black] (7.6,-7.65) circle (1.2pt);
\draw[dotted,line width =0.8pt] (6.7,-7.65)--(7.15,-7.65);
\draw [decorate,line width=0.8] (6.2,-3.8)--(7.6,-3.8);
\node at (5.5,-4.3) {$e$};
\node at (6.9,-3.5) {$k$};
\end{tikzpicture}
\end{center}
\end{figure}
\begin{figure}[hb]
\begin{center}
\begin{tikzpicture}[decoration=brace]
\draw[line width=0.8,pattern=north west lines] (-4,0) arc (0:360:1);
\draw[line width=0.8] (-1,0) arc (0:360:1);
\draw[line width=0.8] (2,0) arc (0:360:1);
\draw[dotted,line width =0.8pt] (2.7,0)--(3.15,0);
\draw[line width=0.8] (6,0) arc (0:360:1);
%\draw[line width =0.8pt] (-2.9,0.4)--(-1.1,0.4);
%\draw[line width =0.8pt] (0.1,0.4)--(1.9,0.4);
%\draw[line width =0.8pt] (4.1,0.4)--(5.9,0.4);
\filldraw[line width=0.8,pattern=north west lines](-2,0)--(-1,0) arc (0:50:1)--(-2,0);
\filldraw[line width=0.8,pattern=north west lines](1,0)--(2,0) arc (0:50:1)--(1,0);
\filldraw[line width=0.8,pattern=north west lines](5,0)--(6,0) arc (0:50:1)--(5,0);
\node at (0,-1.3) {$\mathcal{P}_e$};
\node at (1.5,1.6) {$k$};
\node at (-0.6,0.5) {$\frac{1}{k}$};
\node at (2.4,0.5) {$\frac{1}{k}$};
\node at (6.4,0.5) {$\frac{1}{k}$};
\draw [decorate,line width=0.8] (-3,1.2)--(6,1.2);
\end{tikzpicture}
\end{center}
\caption{The proof sketch of Lemma \ref{cover}}
\end{figure}
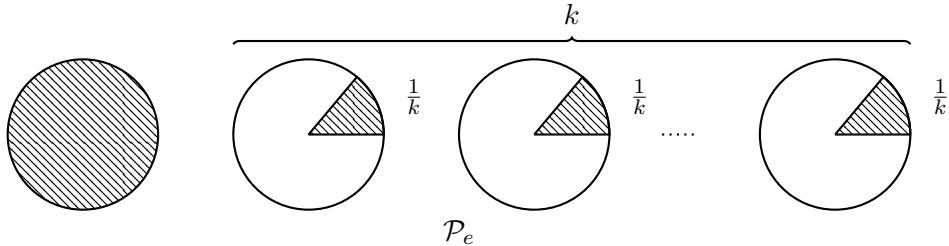

\section{Concluding Remarks}
%The threshold for the minimum $(k-1)$-degree condition in Theorem~\ref{main} is asymptotically equal to the single host hypergraph version in Theorem~\ref{main1}.
Inspired by a series of recent successes on transversal settings of matchings~\cite{2020A,LU,Hongliang2018ON,2020C} and Hamilton cycles~\cite{MR4171383} , we suspect the threshold for transversal Hamilton cycle in a $k$-graph system is the same with the threshold for Hamilton cycle in a single $k$-graph.

%For $k=3$, and large enough hypergraph, R\"{o}dl, Ruci\'{n}ski and Szemer\'{e}di obtained an exact result \cite{EXACT}.
%It is  meaningful to give the exact rainbow version for $3$-graphs.
%Giving the exact co-degree condition for Hamilton cycle in $k$-graph is still open, it causes much more difficulty in the rainbow version.
\begin{conjecture}
Suppose $\textbf{H}=\{H_i\}_{i\in[n]}$ is an $n$-vertex $k$-graph system on $V$, $n\geq k+1\geq 4$, such that $\delta_{k-1}(H_i)\geq \lfloor (n-k+3)/2\rfloor$, then there is an $\textbf{H}$-transversal Hamilton cycle.
\end{conjecture}

On the other hand, the problem of giving the sufficient condition for the transversal Hamilton $\ell$-cycles, $\ell\in[k-2]$, is still open.
\section{Acknowledgement}

This work was supported by the Natural Science Foundation of China (11871311, 11631014) and Youth Interdisciplinary Innovation Group of Shandong University.
We thank the anonymous referee for detailed feedback that improved the presentation of the paper.

%\bibliographystyle{plain}
%\bibliography{1112}

\begin{appendix}
\section{The postponed proofs}
In this section, we prove Lemma~\ref{connecting}.
The idea of the proof is to grow tree-like structures (called cascades) from both designated ends $e_1$ and $e_2$ until they meet, forming the $\emph{\textbf{H}}$-transversal path as desired.
Our proof follows almost identical as that in~\cite{cascade,approximate}.
Before we describe the cascades, it is convenient to introduce the following notation.
For two sequences of vertices
\[
\omega_1=(v_1,\ldots,v_r,w_1,\ldots,w_s)  {\ \rm and\ } \omega_2=(w_1,\ldots,w_s,u_1,\ldots,u_t)
\]
where $r,t\geq1,s\geq0$ and all vertices are distinct, we define their $concatenation$ as
\[
\omega_1\omega_2=(v_1,\ldots,v_r,w_1,\ldots,w_s,u_1,\ldots,u_t).
\]
This operation can be iterated.
For instance, if $\omega_1=(w_1,\ldots,w_{k-2})$, $\omega_2=(w_2,\ldots,w_{k-1})$ and $\omega_3=(w_3,\ldots,w_{k})$ where all $w_i$ are distinct, then $\omega_1\omega_2\omega_3=(w_1,\ldots,w_{k})$.
We could write $\omega_1\omega_2w_k$ instead of $\omega_1\omega_2\omega_3$.
Let $e_0=(v_1,\ldots,v_{k-1})$ be a given $(k-1)$-tuple of vertices.
We will define the \emph{transversal} $e_0$-$cascade$ as an auxiliary sequence of bipartite graphs $G_j, j=1,2,\ldots$, with bipartitions $(A_{j-1},A_j)$, whose vertices are $(k-2)$-tuples of the vertices of $H$ and the edges correspond to some $(k-1)$-tuples of the vertices of $H$.
Each node $f\in A_j$ belongs to two graphs $G_j$ and $G_{j+1}$.
Its neighbors in $G_j$ belongs to $A_{j-1}$, while its neighbors in $G_{j+1}$ belongs to $A_{j+1}$.
For a node $f=(v_1,\ldots,v_{k-2})$ of the transversal cascade, the vertex $v_1$ is called the $prefix$, while $v_{k-2}$ is called the $suffix$ of $f$.

We define the transversal cascade recursively as follows.
Let $e_0=(v_1,\ldots,v_{k-1})$, $f_0=(v_2,\ldots,v_{k-1})$ and $A_0=\{f_0\}$.
For every vertex $v\notin e_0$, we include the node $g=(v_3,\ldots,v_{k-1},v)$ in the set $A_1$ if and only if $v_1f_0g=e_0v\in H_{c_1}$ for $c_1\in[c]$.
The graph $G_1$ is the star with center $f_0$ and the arms leading to all the nodes $g\in A_1$.

Further, let $A_2$ be the set of all $(k-2)$-tuples $h$ such that for some node $g\in A_1$ we have $f_0gh\in H_{c_2}$ where $c_2\neq c_1$ and $c_2\in[c]$.
Note that each $h\in A_2$ is obtained from a node $g\in A_1$ by dropping the prefix of $g$ and adding a new suffix $u$, we denote such node by $g_u$.
The graph $G_2$ consists of all edges $\{g,h\}$ where $g\in A_1$, $h\in A_2$ and $f_0gh\in H_{c_2}$, it is equal to say $G_2$ consists of all edges $\{g,g_u\}$ where $f_0gu\in H_{c_2}$.

For $j=3,\ldots,k-2$, we similarly define
\[
A_j=\{h:\exists f\in A_{j-2}, g\in A_{j-1} {\ \rm such\ that}\ \{f,g\}\in G_{j-1}, fgh\in H_{c_j}\ {\rm where}  \ c_j\neq c_{\ell}\ {\rm for}\ \ell\in[j-1]\}
\]
and $G_j$ as the bipartite graph with bipartition $(A_{j-1}, A_j)$ and the edge set
\[
\{\{g,h\}:\exists f\in A_{j-2}\ {\rm such\ that}\ \{f,g\}\in G_{j-1}\ {\rm and}\ fgh\in H_{c_j}, {\rm \ where}  \ c_j\neq c_{\ell}\ {\rm for}\ \ell\in[j-1]\}.
\]
In other words, $A_j$ and $G_j$ correspond to the sets of $(k-2)$-tuples and $(k-1)$-tuples of the vertices of $V$ which can be reached from $e_0$ in $j$ steps by an $\emph{\textbf{H}}$-transversal path.

\textbf{First refinement.} Having defined $A_j$ and $G_j$ for $j\leq k$, beginning with $j=k-1$ we change the recursive mechanism by getting rid of the nodes in $A_j$ with too small degree in $G_j$.
We define auxiliary
\[
A_{k-1}'=\{h:\exists f\in A_{k-3}, g\in A_{k-2} {\ \rm such\ that}\ \{f,g\}\in G_{k-2}, fgh\in H_{c_{k-1}}\ {\rm where}  \ c_{k-1}\neq c_{\ell}\ {\rm for}\ \ell\in[k-2]\}
\]
and $G_{k-1}'$ as the bipartite graph with bipartition $(A_{k-2}, A_{k-1}')$
and the edge set
\[
\{\{g,h\}:\exists f\in A_{k-3}\ {\rm such\ that}\ \{f,g\}\in G_{k-2}\ {\rm and}\ fgh\in H_{c_{k-1}}\ {\rm where}  \ c_{k-1}\neq c_{\ell}\ {\rm for}\ \ell\in[k-2]\}.
\]
Then let $A_{k-1}$ be the subset of $A_{k-1}'$ consisting of all nodes $h$ with $\deg_{G_{k-1}'}(h)\geq\sqrt{n}$ and set $G_{k-1}=G_{k-1}'[A_{k-2}\cup A_{k-1}]$.
%For convenience, we set $A_j'=A_j$ and $G_j'=G_j$ for all $j\leq k-2$.

\textbf{Second refinement.} For $j\geq k$, to form an edge $\{g,h\}$ of $G_j$ we will now require not one but many nodes $f\in A_{j-2}$ to fulfil the above definition.

Set $m=\lceil n^{1/4}\rceil$. Having defined $G_{j-1}$, let
$A_j'=\{h:\exists f_1,\ldots,f_m\in A_{j-2}, g\in A_{j-1}\ {\rm such\ that\ for\ all}$ $\ i\in[m], \{f_i,g\}\in G_{j-1}\ and \ f_igh\in H_{c_j}\ {\rm where}  \ c_j\neq c_{\ell}\ {\rm for}\ \ell\in[j-1]\}$
 and let $G_j'$ be the bipartite graph with bipartition $(A_{j-1}, A_j')$ and the edge set
$\{\{g,h\}:\exists f_1,\ldots,f_m\in A_{j-2}\ {\rm such\ that\ for\ all}\ i\in[m], \{f_i,g\}\in G_{j-1}\ and \ f_igh\in H_{c_j}\ {\rm where}  \ c_j\neq c_{\ell}\ {\rm for}\ \ell\in[j-1]\}.$

Finally, let $A_j$ be the subset of $A_j'$ consisting of all nodes $h$ with $\deg_{G_j'}(h)\geq\sqrt{n}$ and let $G_j=G_j'[A_{j-1}\cup A_j]$.
The sequence $(G_j),j=1,2,\ldots,$ will be called the transversal $e_0$-cascade.
\subsection{Properties of the cascade}
\begin{claim}[\cite{approximate}]
\label{avoid}
For every $j\geq k-1$ and every edge $\{g,h\}$ of $G_j$ where $g=(w_1,\ldots,w_{k-2})\in A_{j-1}, h=(w_2,\ldots,w_{k-1})\in A_{j}$ and $(g\cup h)\cap e_0=\emptyset$ and for every set of vertices $W\subset V\setminus(g\cup h\cup e_0)$ such that $j+|W|\leq n^{1/4}$, there is an $\textbf H$-transversal path $P$ of length $j$ which connects $(w_{k-1},\ldots,w_1)$ with $e_0=(v_1,\ldots,v_{k-1})$ and $V(P)\cap W=\emptyset$.
\end{claim}
\textbf{Degrees.} Recall that $G_j'=G_j$ for $j\leq k-2$. For a node $g\in A_j$, we set
\[
d^+(g)=\deg_{G_{j+1}'}(g)\ {\rm and}\ d^-(g)=\deg_{G_j}(g)
\]
for the \emph{forward} and \emph{backward} degree of $g$ in the cascade.
Note that in the definition of $d^+(g)$ we consider the forward degree before some small degree vertices of $A_{j+1}'$ are removed.
The reason is that we have no control over the effects of the removal on individual forward degrees.
On the other hand, for all $f\in A_j$, $\deg_{G_j}(f)=\deg_{G_j'}(f)$, so the backward degree is unaffected unless the node is removed.
It is trivial that $d^-(g),d^+(g)\leq n-k+2$.
Observe that $G_1\cup\cdots\cup G_{k-2}$ is a tree, thus, $d^-(g)=1$ for all $g\in A_j,j=1,\ldots,k-2$.
Recall that for $j\geq k-1$ the graph $G_j$ is obtained from $G_j'$ by removing nodes $g$ with $\deg_{G_j'}(g)<\sqrt{n}$.
Hence our construction guarantees that for all $g\in A_j, j\geq k-1$, we have $d^-(g)\geq\sqrt{n}$.

For all $j\leq k-2$ and all $g\in A_j$,
\begin{equation}\label{666}
d^+(g)\geq\left(\frac{1}{2}+\gamma\right)n,
\end{equation}
since there are at least $(\frac{1}{2}+\gamma)n$ vertices $u$ such that $fgu\in H_{c_{j+1}}$ where $f$ is the neighbor of $g$ in $A_{j-1}$.
Each such vertex $u$ corresponds to a neighbor $g_u$ of $g$ in $A_{j+1}$.

For $j\geq k$, the second refinement affects and no lower bound on $d^+(g)$ is obvious.
However, the lower bound $d^-(g)\geq\sqrt{n}$ introduced by the first refinement maintains.
%\begin{claim}
%For all $j\geq k-1$ and each $g\in A_j$ we have
%\[
%d^+(g)\geq(\frac{1}{2}+\gamma)n-n^{3/4}.
%\]
%\end{claim}

%\begin{claim}
%For all nodes $g$ in the cascade
%\begin{itemize}
 % \item if $d^-(g)\geq\frac{1}{2}n$ then $d^+(g)=n-k+2$,
 % \item $d^-(g)\leq d^+(g)$.
%\end{itemize}
%\end{claim}

\textbf{Growth.} By inequality (\ref{666}), for each $j\in[k-2]$, we have
\begin{equation}\label{growth1}
|G_j|=|A_j|\geq\left(\frac{1}{2}+\gamma\right)^jn^j,
\end{equation}
\begin{equation} \label{growth2}
|G_{k-1}|\geq\left(\frac{1}{2}+\gamma\right)^{k-1}n^{k-1}.
\end{equation}

%\begin{claim}
%For all $k-1\leq j\leq n^{1/3}$ we have
%\[
%|G_j|>(\frac{n}{2})^{k-1}
%\]
%and
%\[
%|A_j|>\frac{1}{2}(\frac{n}{2})^{k-2}
%\]
%\end{claim}

Call a node $f\in A_j$ \emph{small} if $d^-(f)<\frac{1}{2}n$ and denote by $S_j$ the subset of $A_j$ consisting of the small nodes.
Assume for simplicity that $1/\varepsilon^2$ is an integer.
\begin{claim}[\cite{approximate}]
\label{5.6}
There exists an index $j_0$, $k-1\leq j_0\leq k-1+(k-1)/\gamma^2$ such that for all $j\in[j_0,j_0+k-2]$ we have $|S_j|\leq2\gamma n^{k-2}$.
\end{claim}
\begin{claim}[\cite{approximate}]
\label{5.7}
Let
\[
k\gamma^{2^{2-k}}<2^{-k}
\]
and let $j_0$ be as in Claim \ref{5.6}. Then
$|A_{j_0+k-2}\backslash S_{j_0+k-2}|\geq(n-k+2-\gamma^{2^{2-k}}n)^{k-2}$.
\end{claim}
\subsection{The completion of the proof of Lemma \ref{connecting}}
Let $\gamma_0$ satisfy the condition in Claim \ref{5.7}, i.e. $\gamma_0:=\gamma^{2^{2-k}}$ and $k\gamma_0<2^{-k}$.
Given two disjoint $(k-1)$-tuples of vertices $e_1$ and $e_2$, we build the transversal $e_1$-cascade and the transversal $e_2$-cascade, with the sets of nodes denoted by $A_j$ and $B_j$.

Let $j_1=j_0+k-2$, where $j_0$ is the index guaranteed by Claim \ref{5.6} for the transversal $e_1$-cascade.
Then by Claim \ref{5.7}, with sufficiently large $n$, using Bernoulli inequality, we have
\[
|A_{j_1}\backslash S_{j_1}|\geq(n-2\gamma_0n)^{k-2}>(1-2k\gamma_0)n^{k-2}.
\]
On the other hand by inequality (\ref{growth1}) for $j=k-2$, we have $|B_{k-2}|>2^{2-k}n^{k-2}$,
\[
|B_{k-2}\cap(A_{j_1}\backslash S_{j_1})|>(2^{2-k}-2k\gamma_0)n^{k-2}\geq\left(\frac{n}{2}\right)^{k-2}.
\]
Hence, there is a not small node $g=(u_1,\ldots,u_{k-2})\in A_{j_1}$ such that $g\cap(e_1\cup e_2)=\emptyset$ and $g'=(u_{k-2},\ldots,u_1)\in B_{k-2}$.

Let $e_2=(w_1,\ldots,w_{k-1})$, $S=\{u_1,\ldots,u_{k-2},w_{k-1}\}$ and $V_0$ be the set of prefixes $v$ of the neighbors $f\in A_{j_1-1}$ of $g$.
Since $g'=(u_{k-2},\ldots,u_1)\in B_{k-2}$, we obtain that $w_1\cdots w_{k-1}u_{k-2}\cdots u_1$ is an $\emph{\textbf{H}}$-transversal path.
By Fact \ref{neighbor}, we have $|N_{H_{c_{j_1}}}(S)\cap V_0|>\gamma n$, and thus, there is at least one vertex $v_0\notin e_2$ such that $\{v_0,u_1,\ldots,u_{k-2},w_{k-1}\}\in H_{c_{j_1}}$.
%Besides, $g'=(u_{k-2},\ldots,u_1)\in B_{k-2}$, which guarantees that there is partial $\emph{\textbf{H}}$-transversal path $u_{k-2}\cdots u_1w_{k-1}\cdots w_1$.

Let $P_1=e_1\cdots v_0u_1\cdots u_{k-2}$ be an $\emph{\textbf{H}}$-transversal path of length $j_1$ which avoids the vertices of $e_2$.
The existence of $P_1$ follows from Claim \ref{avoid} with $W=e_2$.
The path $P$ obtained from $P_1$ by adding the segment $(w_{k-1},\ldots,w_1)$ and the ``hook-up'' edge $\{v_0,u_1,\ldots,u_{k-2},w_{k-1}\}$, is the $\emph{\textbf{H}}$-transversal path connecting $e_1$ and $e_2$ as desired.

By the bound on $j_0$ established in Claim \ref{5.6} and since $\gamma\leq1/2$,
\[
|V(P)|=j_1+2(k-1)=j_0+3k-4\leq\frac{k-1}{\gamma^2}+4k-5\leq\frac{2k}{\gamma^2}.
\]
\begin{figure}[htb]
\centering
%%Requires \usepackage{graphicx}
\includegraphics[width=14cm]{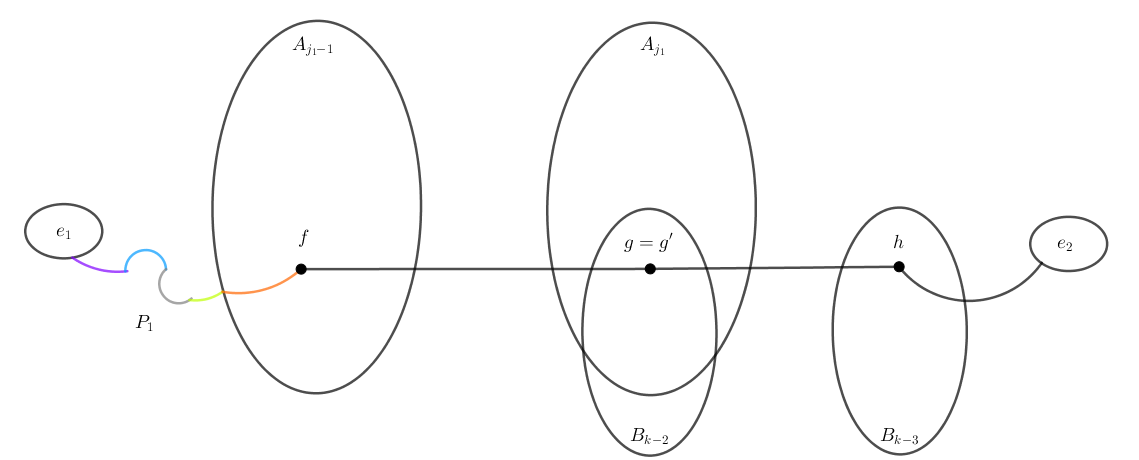}\\
%\caption{Illustration of the rainbow absorbers (with directions omitted)}
%\label{absorber}
\end{figure}
\begin{figure}[htb]
  \centering
 \includegraphics[width=9cm]{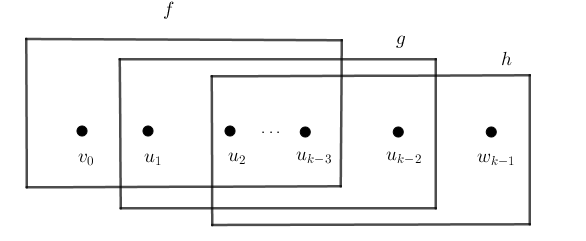}\\
  \caption{An $\emph{\textbf{H}}$-transversal path connecting two $(k-1)$-tuples $e_1$ and $e_2$}%\label{}
\end{figure}

\end{appendix}

\end{document}